\documentclass[12pt]{article}
\usepackage{amssymb, amsmath, amsfonts,dsfont, mathrsfs,euscript,eufrak,colonequals,indentfirst}
\usepackage{textcomp}
\usepackage{latexsym}

\newtheorem{Theorem}{Theorem}

\newtheorem{Proposition}{Proposition}

\newtheorem{Remark}{Remark}

\newtheorem{Assumption}{Assumption}

\newcommand{\N}{\mathbb N}
\newcommand{\R}{\mathbb R}
\newcommand{\Z}{\mathbb Z}

\newcommand{\Expec}{ \mathds{E}\,}

\newcommand{\sign}{\mathrm{sign}}
\newcommand{\dd}{\,d}
\newcommand{\id}{\mathds{1}}

\newcommand{\AppClass}{\mathcal{A}}
\newcommand{\CovFunc}{\mathcal{K}}
\newcommand{\KorElem}{\mathbb{B}}
\newcommand{\e}{\varepsilon}

\newcommand{\mult}{\mathfrak{m}}
\newcommand{\tlambda}{\bar{\lambda}}
\newcommand{\tpsi}{\bar \psi}

\newcommand{\tLambda}{\bar{\Lambda}}
\newcommand{\ContSet}{\mathbf{C}}
\newcommand{\defeq}{\colonequals}
\newcommand{\eqone}{\underset{1}{=}}

\newcommand{\tn}{\bar n}
\newcommand{\CountFunc}{\mathcal{N}}
\newcommand{\Defect}{\mathcal{R}}

\newcommand{\signum}{\mathfrak{s}}

\begin{document}
	
\author{A. A. Khartov\footnote{St. Petersburg State University, 7/9 Universitetskaya nab., 199034 St. Petersburg, Russia, e-mail: \texttt{alexeykhartov@gmail.com}}\qquad M. Zani\footnote{Laboratoire MAPMO, Universit\'e d'Orl\'eans, UFR Sciences, B\^{a}timent Math\'ematiques, Rue de Chartres, B.P. 6759-45067, Orl\'eans cedex 2, France, e-mail: \texttt{marguerite.zani@univ-orleans.fr} }}

\title{Asymptotic analysis of average case approximation complexity of additive random fields}
\maketitle

\begin{abstract}
We study approximation properties of sequences of centered additive random fields $Y_d$, $d\in\N$. The average case approximation complexity $n^{Y_d}(\e)$ is defined as the  minimal number of evaluations of arbitrary linear functionals that is needed to approximate $Y_d$ with relative $2$-average error not exceeding a given threshold $\e\in(0,1)$. We investigate the growth of $n^{Y_d}(\e)$ for arbitrary fixed $\e\in(0,1)$ and $d\to\infty$. Under natural assumptions we obtain general results concerning asymptotics of $n^{Y_d}(\e)$. We apply our results to additive random fields with marginal random processes corresponding to the Korobov kernels.
\end{abstract}


\section{Introduction and problem setting}
Let $X_j(t)$, $t\in[0,1]$, be a given sequence of zero-mean random processes with continuous covariance functions $\CovFunc^{X_j}(t,s)$, $t,s\in[0,1]$, $j\in\N$, where $\N$ is the set of positive integers. For every $d\in\N$ on some probability space we consider a zero-mean random field $Y_d(t)$, $t\in[0,1]^d$, with  the following covariance function:
\begin{eqnarray}\label{def_CovFuncYd}
\CovFunc^{Y_d}(t,s)=\sum_{j=1}^d \CovFunc^{X_j}(t_j,s_j),
\end{eqnarray}
where $t=(t_1,\ldots, t_d)$ and $s=(s_1,\ldots, s_d)$ are from $[0,1]^d$. 
Random fields of such type can be constructed from the marginal processes in the following way. Let the random processes $X_j$, $j\in\N$, be defined on some probability space. Let $X_j$, $j\in\N$, be uncorrelated or independent. Then for every $d\in\N$ the random field
\begin{eqnarray*}
Y_d(t)\defeq\sum_{j=1}^d X_j(t_j),\quad t=(t_1,\ldots, t_d)\in[0,1]^d,
\end{eqnarray*}
will have the covarience function \eqref{def_CovFuncYd}.

Random fields with the described covariance structure belong to wide class of \textit{additive random fields}, whose study is of rather big interest. For example, they are used in approximation of more complicated random fields and appear in the theory of intersections and selfintersections of Brownian processes (see \cite{ChenLi} and the references given there).

We will study approximation properties of defined additive random fields $Y_d$, $d\in\N$. Namely, every  $Y_d(t)$, $t\in[0,1]^d$, is considered as a random element of the space $L_2([0,1]^d)$ with scalar product $(\,\cdot,\,\cdot)_{2,d}$ and  norm $\|\cdot\|_{2,d}$. 
We will investigate the \textit{average case approximation complexity} (simply the \textit{approximation complexity} for short) of $Y_d$, $d\in\N$:
\begin{eqnarray}\label{def_nXde}
n^{Y_d}(\e)\colonequals\min\bigl\{n\in\N:\, e^{Y_d}(n)\leqslant \e\, e^{Y_d}(0)  \bigr \},
\end{eqnarray}
where $\e\in(0,1)$ is a given error threshold, and
\begin{eqnarray*}
e^{Y_d}(n)\colonequals\inf\Bigl\{ \bigl(\Expec\bigl\|Y_d - \widetilde Y^{( n)}_d\bigr\|_{2,d}^2\bigr)^{1/2} :  \widetilde Y^{(n)}_d\in \AppClass_n^{Y_d}\Bigr\}
\end{eqnarray*}
is the smallest 2-average error among all linear approximations of $Y_d$, $d\in\N$, having rank $n\in\N$. The corresponding classes  of linear algorithms are denoted by
\begin{eqnarray*}
	\AppClass_n^{Y_d}\colonequals \Bigl\{\sum_{m=1}^{n} (Y_d,\psi_m)_{2,d}\,\psi_m :  \psi_m \in L_2([0,1]^d)\Bigr \}, \quad d\in\N,\quad n\in\N.
\end{eqnarray*}
We always work with \textit{relative errors}, taking into account the following ``size'' of $Y_d$:
\begin{eqnarray*}
	e^{Y_d}(0)\colonequals  \bigl(\Expec\|Y_d\|_{2,d}^2\bigr)^{1/2}<\infty,
\end{eqnarray*}
which is the approximation error of $Y_d$ by zero element.

The approximation complexity $n^{Y_d}(\e)$ is considered as a function depending on two variables $d\in\N$ and $\e\in(0,1)$. Namely, we will investigate the asymptotic behaviour of $n^{Y_d}(\e)$ for arbirarily small fixed $\e$ and $d\to\infty$.  We assume that covariance characteristics of every marginal process $X_j$, $j\in\N$, are known. More precisely, by assumption, we know eigenvalues and traces of covariance operators of $X_j$, $j\in\N$.

Multivariate approximation problems for additive random fields have been considered in the papers \cite{LifZani1} and \cite{LifZani2} in a variety of settings. However, additive random fields are still less investigated than the tensor product-type ones, whose covariance functions are defined as products of marginal ones (see \cite{LifPapWoz1} and the reference given there). Indeed, in \cite{LifZani1} and \cite{LifZani2} the authors studied only the homogeneous case, where approximated additive random fields constructed (in a special way) from copies of one marginal process. We are not aware any results for the non-homogeneous case, where the random fields are composed of a whole sequence of marginal random processes with generally different covariance functions. For tensor product-type random fields these cases  have been comprehensively studied  within described average case setting in \cite{Khart}. This our work is the first step to this direction for additive random fields.

The paper is organized as follows. In Section 2 we provide necessary preliminaries and formulate main assumptions and basic propositions. In Section 3 we infer an integral representation of the approximation complexity, which will be useful for the next general asymptotic results. The corresponding theorems are formulated and proved in Section 4. Next, in Section 5 for illustration we apply our results to additive random fields with marginal random processes corresponding to the Korobov kernels.

For convenience, throughout the paper we will use the following unified notation for the covariance characteristics of random processes and fields. Let $Z(t)$, $t\in[0,1]^n$, be a given zero-mean random process or field with sample paths from $L_2([0,1]^n)$ with some $n\in\N$. We will denote by $K^Z$ and $\CovFunc^Z$ the covariance operator and the covariance function of $Z(t)$, $t\in[0,1]^n$, respectively. From analytic point of view,  $K^{Z}$ is an integral operator with kernel $\CovFunc^{Z}$, which acts by the following formula
\begin{eqnarray*}
K^{Z} f(t)\defeq \int\limits_{[0,1]^n} \CovFunc^{Z}(t,s) f(s) \dd s,\quad f\in L_2([0,1]^n), \quad t\in[0,1]^n.
\end{eqnarray*}
Let $(\lambda^{Z}_k)_{k\in\N}$ and $(\psi^{Z}_k)_{k\in\N}$ denote the non-increasing sequence of eigenvalues and the corresponding  sequence of eigenvectors of $K^{Z}$ respectively. Then $K^{Z} \psi^{Z}_k(t)=\lambda^{Z}_k\psi^{Z}_k(t)$, $k\in\N$, $t\in[0,1]^n$. Here if $K^Z$ is of rank $p\in\N$, then we formally set $\lambda_k^{Z}\colonequals0$, and $\psi_k^Z\equiv0$ for $k>p$.  Let $\Lambda^{Z}$ denote the trace of $K^{Z}$, i.e. $\Lambda^{Z}\defeq \sum_{k=1}^\infty \lambda^{Z}_k$. 

Also we will use the following notation. We denote by $\R$ the set of real numbers. For any function $f$ we will denote by $\ContSet(f)$ the set of all its continuity points and by $f^{-1}$ the generalized inverse function $f^{-1}(y)\colonequals\inf\bigl\{x\in\R: f(x)\geqslant y\bigr\}$, where $y$ is from the range of $f$. By \textit{distribution function} $F$ we mean a non-decreasing function $F$ on $\R$ that is right-continuous on $\R$, $\lim_{x\to-\infty} F(x)=0$, and $\lim_{x\to\infty} F(x)=1$. We write the expectation of a random variable $Z$ as $\Expec Z$. The relation $a_n\sim b_n$ means that $a_n/b_n\to 1$, $n\to\infty$. The quantity $\id(A)$ equals one for the true logic propositions $A$ and zero for the false ones. The number of elements of a finite set $B$ is denoted by $\#(B)$. The functions $x\mapsto \lfloor x\rfloor$ and $x\mapsto\lceil x\rceil$ are  floor and ceiling functions respectively, i.e. $\lfloor x\rfloor=k\in\Z$  whenever $k\leqslant x< k+1$, and $\lceil x\rceil=k\in\Z$  whenever $k-1< x\leqslant k$. For any numbers $x$ and $y$ the notation $x\eqone y$ means that $y\leqslant x\leqslant y+1$.

\section{Preliminaries and main assumptions}
The approximation complexity  $n^{Y_d}(\e)$ can be described in terms of eigenvalues of $K^{Y_d}$. It is well known (see \cite{WasWoz}) that for any $n\in\N$ the following $n$-rank random field
\begin{eqnarray}\label{def_Xdn}
	\widetilde Y^{(n)}_d(t)\colonequals\sum_{k=1}^n (Y_d,\psi^{Y_d}_k)_{2,d}\, \psi^{Y_d}_k(t),\qquad t\in[0,1]^d,
\end{eqnarray}
minimizes the 2-average case error. Hence formula \eqref{def_nXde} is reduced to
\begin{eqnarray*}
	n^{Y_d}(\e)=\min\Bigl\{n\in\N:\, \Expec 
	\bigl\|Y_d-\widetilde{Y}^{(n)}_d\bigr\|_{2,d}^2\leqslant\e^2\,\Expec \|Y_d\|_{2,d}^2 \Bigr\},\quad  d\in\N,\,\, \e\in(0,1).
\end{eqnarray*}
On account of  \eqref{def_Xdn} and $\Expec (Y_d,\psi^{Y_d}_k)_{2,d}^2=(\psi^{Y_d}_k, K^{Y_d}\psi^{Y_d}_k)_{2,d}=\lambda^{Y_d}_k$, $k\in\N$,  we infer the following representation of the approximation complexity $n^{Y_d}(\e)$ in terms of $\lambda^{Y_d}_k$, $k\in\N$:
\begin{eqnarray}
n^{Y_d}(\e)&=&\min\Bigl\{n\in\N:\, \sum_{k=n+1}^\infty \lambda^{Y_d}_k\leqslant\e^2\,\Lambda^{Y_d} \Bigr\}\nonumber\\
&=&\min\Bigl\{n\in\N: \sum_{k=1}^{n} \lambda^{Y_d}_k\geqslant (1-\e^2)\Lambda^{Y_d}\Bigr\},\quad  d\in\N,\,\, \e\in(0,1).\label{conc_nYde}
\end{eqnarray}
Thus the behaviour of distributions of eigenvalues of $K^{Y_d}$ fully determines the growth of $n^{Y_d}(\e)$ as $d\to\infty$. 

Under the additive structure \eqref{def_CovFuncYd}, the numbers $\lambda^{Y_d}_k$,  $k\in\N$, are generally unknown or not easily depend on $\lambda^{X_j}_k$,  $k\in\N$, $j\in\N$. However, under the following assumption, we can explicitly describe the eigenvalues $\lambda^{Y_d}_k$,  $k\in\N$.\\

\noindent\textbf{Basic assumption.} \textit{For every $j\in\N$ there exist  $\psi_0\in\{\psi^{X_j}_1,\psi^{X_j}_2,\ldots,\psi^{X_j}_k,\ldots\}$ such that $\psi_0(t)=1$ for  all $t\in[0,1]$.}\\

There exist important processes with an eigenvector, which is identically $1$ (see \cite{Hick}). We will make this assumption throughout the paper without saying. 

We now decribe the structure of $(\lambda^{Y_d}_k)_{k\in\N}$, $d\in\N$. For every $j\in\N$ let us denote by $\tlambda^{X_j}_0$ the eigenvalue, which is corresponded to identical $1$. Let $(\tlambda^{X_j}_k)_{k\in\N}$ and $(\tpsi^{X_j}_k)_{k\in\N}$ denote the non-increasing sequence of remaining eigenvalues and the corresponding  sequence of eigenvectors of $K^{X_j}$, respectively. We set $\tLambda^{X_j}\defeq \sum_{k\in\N} \tlambda^{X_j}_k=\Lambda^{X_j}-\tlambda^{X_j}_0$, $j\in\N$. The basic assumption ensures that the family 
\begin{eqnarray}\label{def_psiXjtjk_orthog}
\{1\}\cup\bigl\{\tpsi^{X_j}_k(t_j): k\in\N,\, j=1,\ldots, d,   \bigr\},\quad (t_1,\ldots, t_d)\in[0,1]^d,
\end{eqnarray}
is an orthogonal system in $L_2([0,1]^d)$ for every $d\in\N$. Indeed, it is easily seen that
\begin{eqnarray*}
1, \tpsi^{X_j}_1(t_j), \tpsi^{X_j}_2(t_j),\ldots,\tpsi^{X_j}_k(t_j),\ldots,
\end{eqnarray*}
are orthogonal in $L_2([0,1]^d)$  for every $j=1,\ldots, d$, and $d\in\N$. Next, from the basic assumption we conclude that
\begin{eqnarray*}
\int\limits_{[0,1]} \tpsi^{X_j}_k(t_j)\, \dd t_j =0,\quad k\in \N,\quad j=1,\ldots, d,\quad d\in\N. 
\end{eqnarray*}
Therefore for all $d\in\N$ and $j,l=1,\ldots, d$,  $j\ne l$, we have
\begin{eqnarray*}
\int\limits_{[0,1]^d} \tpsi^{X_j}_k(t_j)\, \tpsi^{X_l}_m(t_l)\dd t=\int\limits_{[0,1]} \tpsi^{X_j}_k(t_j)\, \dd t_j \int\limits_{[0,1]}  \tpsi^{X_l}_m(t_l)\dd t_l=0,\quad k,m \in \N.
\end{eqnarray*}
Thus orthogonality of \eqref{def_psiXjtjk_orthog} is shown. Using this it is easy to check  that for every $d\in\N$ identical 1 is an eigenvector of $K^{Y_d}$ with the eigenvalue $\tlambda^{Y_d}_0\defeq \sum\nolimits_{j=1}^d \tlambda^{X_j}_0$, and that the pairs $\tlambda^{X_j}_m$ and $\tpsi^{X_j}_k(t_j)$, $t_j\in[0,1]$, for all  $j=1,\ldots, d$, and $k\in\N$, are remaining eigenpairs of $K^{Y_d}$. 

For every $d\in\N$ we will denote analogously by $(\tlambda^{Y_d}_k)_{k\in\N}$  the non-increasing sequence of eigenvalues, which are corresponded to eigenvectors of $K^{Y_d}$ from the family \eqref{def_psiXjtjk_orthog} without identical 1.  We also set $\tLambda^{Y_d}\defeq \sum_{m\in\N} \tlambda^{Y_d}_m=\Lambda^{Y_d}-\tlambda^{Y_d}_0$, $d\in\N$. It is easily seen that
\begin{eqnarray*}
\Lambda^{Y_d}=\sum_{j=1}^d\Lambda^{X_j},\qquad \tLambda^{Y_d}=\sum_{j=1}^d \tLambda^{X_j},\quad d\in\N.
\end{eqnarray*}

Let us consider the sequence $\gamma_d=\tlambda^{Y_d}_0/\Lambda^{Y_d}$, $d\in\N$. We exclude the trivial case when $\gamma_d=1$ for every $d\in\N$. If $\gamma_d$ tends to $1$ as $d\to\infty$, then from \eqref{conc_nYde}  we conclude $n^{Y_d}(\e)=1$ for  every $\e\in(0,1)$ and all sufficiently large $d\in\N$. If $\gamma_d$ has not a limit, then the behaviour of $n^{Y_d}(\e)$ can be rather unregular. In order to exclude these cases we will make the following assumption.

\begin{Assumption}\label{assum_Regular}
The sequence 
\begin{eqnarray*}
\dfrac{\tlambda^{Y_d}_0}{\Lambda^{Y_d}}=\dfrac{\sum_{j=1}^d\tlambda^{X_j}_0}{\sum_{j=1}^d\Lambda^{X_j}}, \qquad d\in\N,
\end{eqnarray*}
has a limit, which is strictly less than $1$. 
\end{Assumption}

Let $\e_0$ be the number from $(0,1]$ such that $1-\e_0^2 =\lim\limits_{d\to\infty} \bigl(\tlambda^{Y_d}_0/\Lambda^{Y_d}\bigr)$, i.e.
\begin{eqnarray}\label{def_e0}
1-\e_0^2=\lim_{d\to\infty}\dfrac{\sum_{j=1}^d\tlambda^{X_j}_0}{\sum_{j=1}^d\Lambda^{X_j}},
\end{eqnarray}
or, in other form,
\begin{eqnarray*}
\e_0^2=\lim_{d\to\infty}\dfrac{\sum_{j=1}^d\tLambda^{X_j}}{\sum_{j=1}^d\Lambda^{X_j}}.
\end{eqnarray*}
If $\e_0\in(0,1)$, then for every  $\e\in(\e_0,1)$ the inequality $\tlambda^{Y_d}_0/\Lambda^{Y_d}>1-\e^2$ is satisfied for all sufficiently large $d\in\N$. By the formula \eqref{conc_nYde}, we have  $n^{Y_d}(\e)=1$ for such $\e$ and $d$. In view of these remarks, we will consider the approximation complexity $n^{Y_d}(\e)$ only for $\e\in(0,\e_0)$ under Assumption \ref{assum_Regular}.

\begin{Proposition}\label{pr_nde_toinfty} 
Let $X_j$, $j\in\N$, satisfy Assumption $\ref{assum_Regular}$. 
Let $\e_0\in(0,1]$ be the number from $\eqref{def_e0}$. The following conditions are equivalent:
\begin{itemize}
\item[$(i)$]\quad $\lim\limits_{d\to\infty}n^{Y_d}(\e)=\infty$\quad for every\quad $\e\in(0,\e_0);$
\item[$(ii)$]\quad $\max\limits_{j=1,\ldots,d}\tlambda^{X_j}_1=o\Bigl(\sum_{j=1}^d \tLambda^{X_j}\Bigr)$,\quad $d\to\infty$.
	\end{itemize}
\end{Proposition}
\noindent\textbf{Proof.}\quad $(i)\Rightarrow(ii)$. Suppose that contrary to our claim, there exists a sequence $(d_l)_{l\in\N}$ of positive integers such that 
\begin{eqnarray*}
\max\limits_{j=1,\ldots,d_l}\tlambda^{X_j}_1\sim c\sum_{j=1}^{d_l} \tLambda^{X_j},\quad l\to\infty, \quad\text{for some } c\in(0,1\,].	
\end{eqnarray*}
Since
\begin{eqnarray}\label{conc_tlambdaYd1_tLambdaYd1}
\max\limits_{j=1,\ldots,d}\tlambda^{X_j}_1=\tlambda^{Y_{d}}_1,\quad\text{and}\quad \sum_{j=1}^{d} \tLambda^{X_j}=\Lambda^{Y_{d}}-\tlambda^{Y_{d}}_0,\quad \text{for all}\quad d\in\N,
\end{eqnarray}
we have $\tlambda^{Y_{d_l}}_1\sim c\,\e_0^2 \Lambda^{Y_{d_l}}$ and, consequently, $\tlambda^{Y_{d_l}}_0+\tlambda^{Y_{d_l}}_1\sim\bigl(1-(1-c)\,\e_0^2\bigr)\Lambda^{Y_{d_l}}$ as $l\to\infty$, by Assumption 1. Next, we choose $\e\in (0,1)$ to satisfy $(1-c)\,\e_0^2<\e^2<\e_0^2$. Then $\tlambda^{Y_{d_l}}_0+\tlambda^{Y_{d_l}}_1> (1-\e^2)\Lambda^{Y_{d_l}}$ for all sufficiently large $l$. This gives $n^{Y_{d_l}}(\e)\leqslant 2$, which contradicts $(i)$.

$(i)\Leftarrow(ii)$. Let us arbitrarily fix $\e\in(0,\e_0)$, where $\e_0\in(0,1]$. On account of Assumption~\ref{assum_Regular}, $(ii)$ and \eqref{conc_tlambdaYd1_tLambdaYd1}, we have $\tlambda^{Y_d}_0\sim (1-\e_0^2)\Lambda^{Y_d}$ and $\tlambda^{Y_d}_1=o(\Lambda^{Y_d})$, $d\to\infty$. If $\e_0=1$, then $\tlambda^{Y_d}_0=o(\Lambda^{Y_d})$, $d\to\infty$, and $(i)$ follows from the formula \eqref{conc_nYde} and the inequality
\begin{eqnarray*}
n^{Y_d}(\e)\geqslant \sum_{m=1}^{n^{Y_d}(\e)}\dfrac{\lambda^{Y_d}_m}{\lambda^{Y_d}_1}=\sum_{m=1}^{n^{Y_d}(\e)}\dfrac{\lambda^{Y_d}_m}{\max\{\tlambda^{Y_d}_0,\tlambda^{Y_d}_1\}}\geqslant \dfrac{(1-\e^2)\Lambda^{Y_d}}{\max\{\tlambda^{Y_d}_0,\tlambda^{Y_d}_1\}}.
\end{eqnarray*}
If $\e_0<1$, then for all sufficiently large $d$ the eigenvalue $\tlambda^{Y_d}_0$ is the largest  of eigenvalues $\lambda^{Y_d}_k$, $k\in\N$, i.e. $\lambda^{Y_d}_1=\tlambda^{Y_d}_0$ and, consequently, $\lambda^{Y_d}_2=\tlambda^{Y_d}_1$.  By Assumption \ref{assum_Regular}, for all sufficiently large $d\in\N$ we have $\tlambda^{Y_d}_0<(1-\e^2)\Lambda^{Y_d}$, i.e. $n^{Y_d}(\e)\geqslant 2$. Therefore the following inequalities hold
\begin{eqnarray*}
	n^{Y_d}(\e)-1\geqslant \sum_{m=2}^{n^{Y_d}(\e)}\dfrac{\lambda^{Y_d}_m}{\lambda^{Y_d}_2}
	\geqslant\dfrac{(1-\e^2)\Lambda^{Y_d}-\lambda^{Y_d}_1}{\lambda^{Y_d}_2}
	=\dfrac{(1-\e^2)\Lambda^{Y_d}-\tlambda^{Y_d}_0}{\tlambda^{Y_d}_1}.
\end{eqnarray*}
The last quantity is equivalent to $(\e_0^2-\e^2)\Lambda^{Y_d}/\tlambda^{Y_d}_1$ as $d\to\infty$, where the latter tends to infinity. This implies $(i)$. \quad $\Box$\\

On account of Proposition \ref{pr_nde_toinfty}, in order to obtain the asymptotic results concerning $n^{Y_d}(\e)$, we have to make the following asumption, which seems, however, rather weak.

\begin{Assumption}\label{assum_tendstoinfty}
The following equality holds
\begin{eqnarray*}
\max\limits_{j=1,\ldots,d}\tlambda^{X_j}_1=o\Bigl(\sum_{j=1}^d \tLambda^{X_j}\Bigr),\quad d\to\infty.
\end{eqnarray*}
\end{Assumption}

\begin{Remark}\label{rem_sumLambda1diverg}
Under Assumption $\ref{assum_tendstoinfty}$ the series $\sum_{j=1}^\infty \tLambda^{X_j}$  diverges.
\end{Remark}
\textbf{Proof.}\quad Suppose, contrary to our claim, that $\sum_{j=1}^\infty \tLambda^{X_j}$ converges. Set $d_0\defeq \min \{j\in\N:\lambda^{X_j}_1>0\}$. There exists $d_1>d_0$ such that  $\tLambda^{X_d}<\tlambda^{X_{d_0}}_1$ for all $d\geqslant d_1$. This gives
\begin{eqnarray*}
	\max\limits_{j=1,\ldots,d}\tlambda^{X_j}_1=\max\limits_{j=1,\ldots,d_1}\tlambda^{X_j}_1\quad  \text{for all}\quad d\geqslant d_1.
\end{eqnarray*}
The latter maximum is not smaller than $\lambda^{X_{d_0}}_1>0$. Hence equality in Assumption \ref{assum_tendstoinfty} could not be satisfied, a contradiction.\quad $\Box$

Let us formulate a useful sufficient condition of boundedness of the approximation complexity $n^{Y_d}(\e)$ on $d\in\N$ for every fixed $\e\in(0,1)$.

\begin{Proposition}\label{pr_nYdbounded}
If the series $\sum_{j=1}^\infty \tLambda^{X_j}$  converges, then
\begin{eqnarray}\label{pr_nYdbounded_conc}
\sup_{d\in\N} n^{Y_d}(\e)<\infty,\quad \text{for all}\quad \e\in(0,1).
\end{eqnarray}
\end{Proposition}
\textbf{Proof.}\quad Let us consider the series $\sum_{j=1}^\infty\tlambda^{X_j}_0$. If this diverges, then the sequence
\begin{eqnarray*}
\dfrac{\tlambda^{Y_d}_0}{\Lambda^{Y_d}}=\dfrac{\sum_{j=1}^d\tlambda^{X_j}_0}{\sum_{j=1}^d\tlambda^{X_j}_0+\sum_{j=1}^d\tLambda^{X_j}},\quad d\in\N,
\end{eqnarray*}
tends to $1$. By the remarks before Assumption \ref{assum_Regular}, we have $n^{Y_d}(\e)=1$ for all sufficiently large $d\in\N$. This implies \eqref{pr_nYdbounded_conc}. If the series $\sum_{j=1}^\infty\tlambda^{X_j}_0$ converges, then Assumption~\ref{assum_Regular} holds. Let $\e_0\in(0,1]$ be the number from \eqref{def_e0}. Next, we find $d_0\in\N$ such that
\begin{eqnarray}\label{conc_tLambdaYd}
0<\tLambda^{Y_d}/\Lambda^{Y_d}=1-\tlambda^{Y_d}_0/\Lambda^{Y_d}\leqslant 2\e^2_0,\quad\text{for all}\quad d\in\N,\quad d\geqslant d_0.
\end{eqnarray}

Let $\tlambda^{\infty}_m$, $m\in\N$, be non-increasing  sequence, which consists of all numbers $\tlambda^{X_j}_k$,  $j\in\N$, $k\in\N$. On account of the structure of $\tlambda^{Y_d}_k$, $k\in\N$, $d\in\N$, we have 
\begin{eqnarray*}
\dfrac{1}{\tLambda^{Y_d}} \sum\limits_{k=N}^{\infty}\tlambda^{Y_d}_k\leqslant \dfrac{1}{\tLambda^{Y_{d_0}}} \sum\limits_{k=N}^{\infty}\tlambda^{Y_d}_k\leqslant \dfrac{1}{\tLambda^{Y_{d_0}}} \sum\limits_{k=N}^{\infty}\tlambda^{\infty}_k,\quad\text{for all}\quad d\in\N,\quad d\geqslant d_0. 
\end{eqnarray*}
Since $\sum_{j=1}^\infty \tLambda^{X_j}$  converges and $\tlambda^{\infty}_k$, $k\in\N$, are non-negative, the series  $\sum\limits_{k=1}^{\infty}\tlambda^{Y_\infty}_k$  also converges. From this and the previous inequalities we conclude that
\begin{eqnarray*}
\sup_{d\in\N, d\geqslant d_0}\dfrac{1}{\tLambda^{Y_d}} \sum\limits_{k=N}^{\infty}\tlambda^{Y_d}_k\to 0,\quad N\to\infty.
\end{eqnarray*}
Let us choose $N_\e\in\N$ such that 
\begin{eqnarray}\label{conc_Ne}
\sum\limits_{k=1}^{N_\e}\tlambda^{Y_d}_k\geqslant \bigl(1-(\e/\e_0)^2/2\bigr)\tLambda^{Y_d},\quad\text{for all}\quad d\in\N,\quad d\geqslant d_0.
\end{eqnarray}
For these $d$ the inequalities \eqref{conc_tLambdaYd} and \eqref{conc_Ne} imply 
\begin{eqnarray*}
\tlambda^{Y_d}_0+\sum\limits_{k=1}^{N_\e}\tlambda^{Y_d}_k&=&\Lambda^{Y_d}-\tLambda^{Y_d}+\sum\limits_{k=1}^{N_\e}\tlambda^{Y_d}_k\\
&\geqslant& \Lambda^{Y_d}-\tLambda^{Y_d}+\bigl(1-(\e/\e_0)^2/2\bigr)\tLambda^{Y_d}\\
&=& \Lambda^{Y_d}-(\e/\e_0)^2\tLambda^{Y_d}/2\\
&\geqslant& (1-\e^2)\Lambda^{Y_d}.
\end{eqnarray*}
By the formula \eqref{conc_nYde} for $n^{Y_d}(\e)$, we have $n^{Y_d}(\e)\leqslant 1+N_\e$ for all $d\in\N$, $d\geqslant d_0$. Obviously, for every $\e\in(0,1)$ there exists a number $M_\e>0$ such that $n^{Y_d}(\e)<M_\e$ for all $d=1,\ldots, d_0$. Thus \eqref{pr_nYdbounded} is proved.\quad $\Box$

\section{Integral representation of approximation complexity}
Here we will infer a useful integral representation of $n^{Y_d}(\e)$. This will be a base for  our next results.

Let $(\e_d)_{d\in\N}$ be the sequence of numbers from $[0,1]$ such that 
\begin{eqnarray}\label{def_ed}
\e_d^2=\dfrac{\sum_{j=1}^d\tLambda^{X_j}}{\sum_{j=1}^d\Lambda^{X_j}}=1-\dfrac{\sum_{j=1}^d\tlambda^{X_j}_0}{\sum_{j=1}^d\Lambda^{X_j}},\quad d\in\N.
\end{eqnarray}
As we said before Assumption \ref{assum_Regular}, we always suppose that $\e_d>0$ for all $d\in\N$ perhaps except a finite number. Let us introduce the following important sequence of distribution functions:
\begin{eqnarray}\label{def_Fd}
F_d(x)\defeq\dfrac{\sum_{j=1}^d\sum_{k=1}^\infty \tlambda^{X_j}_k \id\bigl(\tlambda^{X_j}_k\geqslant e^{-x}\bigr) }{\sum_{j=1}^d \tLambda^{X_j}},\qquad x\in\R,\quad d\in\N.
\end{eqnarray}
The next theorem shows that for every $d\in\N$ the approximation complexity $n^{Y_d}(\e)$ is in fact fully determined by the function $F_d$. Moreover, we will see below that the sequence $F_d$, $d\in\N$, is convenient to obtain asymptotics of $n^{Y_d}(\e)$ as $d\to\infty$.

\begin{Theorem}\label{th_nYd_repr}
For every $d\in\N$ such that $\e_d>0$, and for every $\e\in(0,\e_d)$ the approximation complexity $n^{Y_d}(\e)$ admits the following representation
\begin{eqnarray}\label{th_nYd_repr_conc}
n^{Y_d}(\e)\eqone\Biggl\lceil \sum_{j=1}^d \tLambda^{X_j}\cdot\!\!\! \int\limits_{0}^{1-(\e/\e_d)^2} \exp\bigl\{F^{-1}_d(y)\bigr\}\dd y \Biggr\rceil.
\end{eqnarray}
\end{Theorem}
\textbf{Proof.}\quad Fix $d\in\N$ such that $\e_d>0$ and $\e\in(0,\e_d)$. Let us consider the quantity
\begin{eqnarray*}
\tn^{Y_d}(\e)\defeq\min\Bigl\{n\in\N: \tlambda^{Y_d}_0 +\sum_{m=1}^{n} \tlambda^{Y_d}_m\geqslant (1-\e^2)\Lambda^{Y_d}\Bigr\}.
\end{eqnarray*}
If $\tlambda^{Y_d}_0\in\{\lambda^{Y_d}_1,\ldots, \lambda^{Y_d}_{n^{Y_d}(\e)}\}$, then $\tn^{Y_d}(\e)=n^{Y_d}(\e)-1$. If $\tlambda^{Y_d}_0\notin\{\lambda^{Y_d}_1,\ldots, \lambda^{Y_d}_{n^{Y_d}(\e)}\}$, i.e.  $\tlambda^{Y_d}_0<\lambda^{Y_d}_{n^{Y_d}(\e)}$, then $\tn^{Y_d}(\e)= n^{Y_d}(\e)-1$ or $\tn^{Y_d}(\e)= n^{Y_d}(\e)$. Thus $n^{Y_d}(\e)\eqone \tn^{Y_d}(\e)$.

Let us represent $\tn^{Y_d}(\e)$ in the following form:
\begin{eqnarray*}
\tn^{Y_d}(\e)&=&\min\Bigl\{n\in\N: \sum_{m=1}^{n} \tlambda^{Y_d}_m\geqslant (1-\e^2)\Lambda^{Y_d}-\tlambda^{Y_d}_0\Bigr\}\\
&=&\min\Bigl\{n\in\N: \sum_{m=1}^{n} \tlambda^{Y_d}_m\geqslant (\e_d^2-\e^2)\Lambda^{Y_d}\Bigr\}\\
&=&\min\Bigl\{n\in\N: \sum_{m=1}^{n} \tlambda^{Y_d}_m\geqslant \bigl(1-(\e/\e_d)^2\bigr)\tLambda^{Y_d}\Bigr\}.
\end{eqnarray*}
It is not difficult to check that
\begin{eqnarray}\label{conc_nYde_CountFuncDefect}
\tn^{Y_d}(\e)=\bigl\lceil\CountFunc_d(\e)-\Defect_d(\e) \bigr\rceil,
\end{eqnarray}
where
\begin{eqnarray*}
\CountFunc_d(\e)&\defeq& \#\bigl\{k\in\N: \tlambda^{Y_d}_k\geqslant \tlambda^{Y_d}(\e)\bigr\},\\
\Defect_d(\e)&\defeq& \dfrac{1}{\tlambda^{Y_d}(\e)} \Bigl(\sum_{k=1}^{\infty} \tlambda^{Y_d}_k\,\id\bigl(\tlambda^{Y_d}_k\geqslant \tlambda^{Y_d}(\e)\bigr)-\bigl(1-(\e/\e_d)^2\bigr)\tLambda^{Y_d}\Bigr),\\
\tlambda^{Y_d}(\e)&\defeq& \tlambda^{Y_d}_{\tn^{Y_d}(\e)}.
\end{eqnarray*}

First observe that
\begin{eqnarray*}
-\ln\tlambda^{Y_d}(\e)
&=& \inf\Bigl\{x\in\R: \sum_{k=1}^\infty\tlambda^{Y_d}_k\,\id\bigr(\tlambda^{Y_d}_k\geqslant e^{-x}\bigl)\geqslant \bigl(1-(\e/\e_d)^2\bigr)\tLambda^{Y_d}\Bigr\}\\
&=& \inf\Bigl\{x\in\R: \sum_{j=1}^d\sum_{k=1}^\infty\tlambda^{X_j}_k\,\id\bigr(\tlambda^{X_j}_k\geqslant e^{-x}\bigl)\geqslant \bigl(1-(\e/\e_d)^2\bigr)\sum_{j=1}^d\tLambda^{X_j}\Bigr\}\\
&=& \inf\Bigl\{x\in\R: F_d(x)\geqslant 1-(\e/\e_d)^2\Bigr\}\\
&=& F^{-1}_d\bigl(1-(\e/\e_d)^2\bigr).
\end{eqnarray*}
Define $p_d(\e)\defeq F_d\bigl( -\ln\tlambda^{Y_d}(\e)\bigr)$. Since $F_d$ is discontinuous at $(-\ln\tlambda^{Y_d}(\e))$, we have the inequality $F_d\bigl( -\ln\tlambda^{Y_d}(\e)-0\bigr)<p_d(\e)$. By the right-continuity of $F_d$, we have $p_d(\e)\geqslant 1-(\e/\e_d)^2$.
Therefore
\begin{eqnarray}\label{conc_invFd_values1e2}
F^{-1}_d(p)=-\ln\tlambda^{Y_d}(\e),\quad 1-(\e/\e_d)^2\leqslant p\leqslant p_d(\e).
\end{eqnarray}

Let $V_d=\{v_{d,1}, v_{d,2},\ldots,v_{d,k},\ldots\}$ be the set of all positive values of the eigenvalues $\tlambda^{Y_d}_k, k\in\N$. Here we assume that $v_{d,1}> v_{d,2}>\ldots>v_{d,k}>\ldots$\,. The set $V_d$ can be finite or infinite, and in the latter case we formally write $\#(V_d)=\infty$. Next, we define $\mult_{d,k}\defeq\#\{m\in\N:\tlambda^{Y_d}_m=v_{d,k}\}$  and $x_{d,k}\defeq -\ln v_{d,k}$,  $k=1,2,\ldots,\#(V_d)$. Observe that for these $k$  we have
\begin{eqnarray}\label{conc_Fd_jumps}
\mult_{d,k}\,\dfrac{v_{d,k}}{\tLambda^{Y_d}}=F_d(x_{d,k})-F_d(x_{d,k-1}),
\end{eqnarray}
and also
\begin{eqnarray}\label{conc_QYd_values}
x_{d,k}=F^{-1}_d(p),\quad p\in \bigl(F_d(x_{d,k-1}),F_d(x_{d,k})\bigr],
\end{eqnarray}
where we formally set $x_{d,0}\defeq-\infty$ and $F_d(x_{d,0})\defeq 0$.

We now consider $\CountFunc_d(\e)$:
\begin{eqnarray*}
\CountFunc_d(\e)=\!\!\!\sum_{\substack{k\in\N:\\\tlambda^{Y_d}_k\geqslant \tlambda^{Y_d}(\e) }}\!\!\! 1
	=\!\!\!\sum_{\substack{k\in\N:\\\tlambda^{Y_d}_k\geqslant \tlambda^{Y_d}(\e) }}\biggl(\dfrac{\tLambda^{Y_d}}{\tlambda^{Y_d}_k}\cdot\dfrac{\tlambda^{Y_d}_k}{\tLambda^{Y_d}}\biggr).
\end{eqnarray*}
Rewrite this in terms of $v_{d,k}$:
\begin{eqnarray*}
\CountFunc_d(\e)=\sum_{k=1}^{\#(V_d)}\biggl(\dfrac{\tLambda^{Y_d}}{ v_{d,k}}\cdot\mult_{d,k}\,\dfrac{v_{d,k}}{\tLambda^{Y_d}}\biggr)\,\id\bigl(v_{d,k}\geqslant \tlambda^{Y_d}(\e)\bigr).
\end{eqnarray*}
Using \eqref{conc_Fd_jumps} and \eqref{conc_QYd_values} we get
\begin{eqnarray*}
\CountFunc_d(\e)&=&\tLambda^{Y_d}\sum_{k=1}^{\#(V_d)}e^{x_{d,k}}\,\bigl(F_d(x_{d,k})-F_d(x_{d,k-1})\bigr)\,\id\bigl(x_{d,k}\leqslant -\ln\tlambda^{Y_d}(\e)\bigr)\\
&=&\tLambda^{Y_d}\sum_{k=1}^{\#(V_d)}\Biggl(e^{x_{d,k}}\int\limits_{F_d(x_{d,k-1})}^{F_d(x_{d,k})} \,\dd y\cdot\id\bigl(x_{d,k}\leqslant -\ln\tlambda^{Y_d}(\e)\bigr)\Biggr)\\
&=&\tLambda^{Y_d}\sum_{k=1}^{\#(V_d)}\Biggl(\,\,\int\limits_{F_d(x_{d,k-1})}^{F_d(x_{d,k})} e^{F^{-1}_d(y)}\,\dd y\cdot\id\bigl(x_{d,k}\leqslant -\ln\tlambda^{Y_d}(\e)\bigr)\Biggr).
\end{eqnarray*}
Since $ -\ln\tlambda^{Y_d}(\e)\in\bigl\{x_{d,k}: k=1,2,\ldots, \#(V_d)\bigr\}$, we have
\begin{eqnarray}\label{conc_CountFuncde}
\CountFunc_d(\e)=\tLambda^{Y_d}\int\limits_{0}^{p_d(\e)} \exp\bigl\{F^{-1}_d(y)\bigr\}\,\dd y=\sum_{j=1}^d\tLambda^{X_j}\cdot\int\limits_{0}^{p_d(\e)} \exp\bigl\{F^{-1}_d(y)\bigr\}\,\dd y.
\end{eqnarray}

We next consider $\Defect_d(\e)$. According to structure of $\tlambda^{Y_d}_k$, $k\in\N$, and definition of $F_d$, we can write 
\begin{eqnarray*}
\Defect_d(\e)&=&\dfrac{1}{\tlambda^{Y_d}(\e)} \Bigl(\sum_{j=1}^d\sum_{k=1}^{\infty} \tlambda^{X_j}_k\,\id\bigl(\tlambda^{X_j}_k\geqslant \tlambda^{Y_d}(\e)\bigr)-\bigl(1-(\e/\e_d)^2\bigr)\sum_{j=1}^d\tLambda^{X_j}\Bigr),\\
&=&\dfrac{\sum_{j=1}^d\tLambda^{X_j} }{\tlambda^{Y_d}(\e)}\,\Bigl(F_d\bigl(-\ln\tlambda^{Y_d}(\e)\bigr)-\bigl(1-(\e/\e_d)^2\bigr)\Bigr)\\
&=&\dfrac{\sum_{j=1}^d\tLambda^{X_j} }{\tlambda^{Y_d}(\e)}\,\Bigl(p_d(\e)-\bigl(1-(\e/\e_d)^2\bigr)\Bigr)\\
&=&\sum_{j=1}^d\tLambda^{X_j}\cdot \exp\bigl\{-\ln\tlambda^{Y_d}(\e)\}\cdot\!\!\!\int\limits_{1-(\e/\e_d)^2}^{p_d(\e)}\dd y.
\end{eqnarray*}
On account of \eqref{conc_invFd_values1e2}, we have
\begin{eqnarray}\label{conc_Defectde}
\Defect_d(\e)=\sum_{j=1}^d\tLambda^{X_j}\cdot\!\!\! \int\limits_{1-(\e/\e_d)^2}^{p_d(\e)} \exp\bigl\{F^{-1}_d(y)\bigr\}\dd y.
\end{eqnarray}
Substituting \eqref{conc_CountFuncde} and \eqref{conc_Defectde} into \eqref{conc_nYde_CountFuncDefect}, we obtain
\begin{eqnarray*}
\tn^{Y_d}(\e)=\Biggl\lceil \sum_{j=1}^d\tLambda^{X_j}\cdot\!\!\!  \int\limits_{0}^{1-(\e/\e_d)^2} \exp\bigl\{F^{-1}_d(y)\bigr\}\dd y \Biggr\rceil.
\end{eqnarray*}
The required representation of $n^{Y_d}(\e)$ immediately follows from the equality $n^{Y_d}(\e)\eqone \tn^{Y_d}(\e)$, which has already been got above.\quad $\Box$

This theorem yields in fact an exact expression for $n^{Y_d}(\e)$ in the homogeneous case.
Namely, we assume $\CovFunc^{X_1}=\CovFunc^{X_2}=\ldots=\CovFunc^{X_j}=\ldots$ . Therefore $X_j$, $j\in\N$, have the same covariance operators and their eigenvalues. Let us denote $\tlambda_k\defeq \tlambda^{X_j}_k$, $k\in\N\cup\{0\}$, and $\Lambda\defeq \Lambda^{X_j}$, $\tLambda\defeq \tLambda^{X_j}$. Thus we have
\begin{eqnarray*}
\sum_{j=1}^d\tLambda^{X_j}=d \tLambda,\quad d\in\N.
\end{eqnarray*}
Observe that the sequence $(\e_d)_{d\in\N}$ is constant:
\begin{eqnarray*}
\e_d^2=\dfrac{\sum_{j=1}^d\tLambda^{X_j}}{\sum_{j=1}^d\Lambda^{X_j}}=\dfrac{\tLambda}{ \Lambda},\quad d\in\N.
\end{eqnarray*}
Denote $\e_0\defeq (\tLambda/\Lambda)^{1/2}$, i.e. $\e_d=\e_0$, $d\in\N$. It is easily seen that all $F_d$ are equal to the following function
\begin{eqnarray*}
F(x)\defeq\dfrac{1}{\tLambda}\sum_{k=1}^\infty \tlambda_k \id\bigl(\tlambda_k\geqslant e^{-x}\bigr),\quad x\in\R.
\end{eqnarray*}
  
Thus from \eqref{th_nYd_repr_conc}  we infer the following expression of the approximation complexity $n^{Y_d}(\e)$ for the homogeneous case:
\begin{eqnarray*}
n^{Y_d}(\e)\eqone\Biggl\lceil d \tLambda\cdot\!\!\!\!\! \int\limits_{0}^{1-(\e/\e_0)^2} \exp\bigl\{F^{-1}(y)\bigr\}\dd y\Biggr\rceil,\quad \e\in(0,\e_0), \quad d\in\N.
\end{eqnarray*}
We see that $n^{Y_d}(\e)$ grows on $d$ as the linear function $d\mapsto C_\e d$,  where the component $C_\e$ is constant for fixed $\e\in(0,\e_0)$.
 
In the non-homogeneous case the formula \eqref{th_nYd_repr_conc} does not give directly  the character of growth of $n^{Y_d}(\e)$ on $d$, because for every $d\in\N$ the function $F_d$ can be rather difficult. The problem is solved by asymptotic analysis  of $n^{Y_d}(\e)$ as $d\to\infty$, under some regular behaviour of distributions of eigenvalues $\lambda^{X_j}_k$, $k\in\N$, $j\in\N$.

\section{Asymptotic analysis of approximation complexity}
In this section we will obtain general results concerning asymptotics of $n^{Y_d}(\e)$ for fixed $\e$ and $d\to\infty$. Theorems \ref{th_nYd_UWconv} and \ref{th_nYd_UWconv_log} below show that if for large $j\in\N$ the distribution of $\tlambda^{X_j}_k$, $k\in\N$, behaves regularly (say as $U$), and if for large $d\in\N$ the distribution of weights $\tLambda^{X_1},\ldots, \tLambda^{X_d}$ (in the common ``size'' $\tLambda^{Y_d}=\sum_{j=1}^d \tLambda^{X_j}$ of $Y_d$) is also regular (say as $W$), then the growth of the approximation complexity $n^{Y_d}(\e)$ is regular, and we can find its asymptotics as $d\to\infty$. The form of the asymptotics depends on the distribution functions $U$ and $W$.

\begin{Theorem}\label{th_nYd_UWconv}
Let $X_j$, $j\in\N$, satisfy Assumptions $\ref{assum_Regular}$ and $\ref{assum_tendstoinfty}$. 
	Let $\e_0\in(0,1]$ be the number from $\eqref{def_e0}$. Suppose that
	\begin{eqnarray}\label{th_nYd_UWconv_cond_U}
	\dfrac{1  }{\tLambda^{X_j}}\sum_{k=1}^\infty\tlambda^{X_j}_k\,\id\bigl(-\ln\tlambda^{X_j}_k-\ell_j\leqslant x\bigr)\to U(x),\quad j\to\infty,\quad\text{for all}\quad x\in\ContSet(U),
	\end{eqnarray}
	with a distribution function $U$ and a monotonic sequence $(\ell_j)_{j\in\N}$.
	Let $\signum=1$ if $(\ell_j)_{j\in\N}$ is non-decreasing, and  $\signum=-1$ if $(\ell_j)_{j\in\N}$ is non-increasing. Next, suppose that 
	\begin{eqnarray}\label{th_nYd_UWconv_cond_W}
	\dfrac{\sum_{j=1}^d\tLambda^{X_j} \id\bigl(\signum\cdot(\ell_j-a_d)\leqslant x\bigr)}{\sum_{m=1}^d\tLambda^{X_m}}\to W(x),\quad d\to\infty,\quad\text{for all}\quad x\in\ContSet(W),
	\end{eqnarray}
with a distribution function $W$ and a sequence $(a_d)_{d\in\N}$.	
	Then for every $\e\in(0,\e_0)$
	\begin{eqnarray}\label{th_nYd_UWconv_conc}
		n^{Y_d}(\e)\sim e^{a_d}\sum_{j=1}^d \tLambda^{X_j}\!\!\!\!\! \int\limits_{0}^{1-(\e/\e_0)^2} \exp\bigl\{F^{-1}(y)\bigr\}\dd y,\qquad d\to\infty,
	\end{eqnarray}
	where $F(x)=\int\limits_{-\infty}^\infty U(x-\signum\cdot v)\dd W(v)$, $x\in\R$. 
\end{Theorem}
\textbf{Proof.}\quad Let $\e_d$ and $F_d$, $d\in\N$, be defined by \eqref{def_ed} and \eqref{def_Fd} respectively.  By Assumption \ref{assum_Regular}, we have $\e_d\to\e_0$ as $d\to\infty$. 
Fix $\e\in (0,\e_0)$. Thus Theorem~\ref{th_nYd_repr} yields the representation \eqref{th_nYd_repr_conc} of $n^{Y_d}(\e)$ for all sufficiently large $d\in\N$. According to Assumption \ref{assum_tendstoinfty} and Proposition \ref{pr_nde_toinfty}, we have $n^{Y_d}(\e)\to\infty$ as $d\to\infty$. In particular, it follows that
\begin{eqnarray}\label{conc_nYdesimint}
n^{Y_d}(\e)\sim\sum_{j=1}^d \tLambda^{X_j}\!\!\!\!\!\int\limits_{0}^{1-(\e/\e_d)^2} \exp\bigl\{F^{-1}_d(y)\bigr\}\dd y,\quad d\to\infty.
\end{eqnarray}

Let us consider the sequence of distribution functions $F_d$, $d\in\N$. We first show that
\begin{eqnarray}\label{conc_FdtoFad}
F_d(x+a_d)\to F(x),\quad d\to\infty,\quad \text{for all}\quad x\in\ContSet(F),
\end{eqnarray}
with the required $F$. Define the distribution functions
\begin{eqnarray}\label{def_Uj}
U_j(x)\defeq\dfrac{1  }{\tLambda^{X_j}}\sum_{k=1}^\infty\tlambda^{X_j}_k\,\id\bigl(\tlambda^{X_j}_k\geqslant e^{-x}\bigr),\quad x\in\R,\quad j\in\N.
\end{eqnarray}
Let us represent $F_d$, $d\in\N$, in the following form
\begin{eqnarray*}
	F_d(x)=\sum_{j=1}^d\Biggl(\dfrac{\tLambda^{X_j} }{\sum_{m=1}^d\tLambda^{X_m}}\cdot U_j(x)\Biggr),\quad x\in\R,\quad d\in\N.
\end{eqnarray*}
For every $d\in\N$ we define
\begin{eqnarray}\label{def_wd_nud}
	w_d(y)\defeq\sum_{k=1}^d\dfrac{\tLambda^{X_k} \id\bigl(k\leqslant y\bigr)}{\sum_{m=1}^d\tLambda^{X_m}},\quad y\in\R,\qquad\text{and}\qquad \nu_d(z)\defeq w^{-1}_d(z),\quad z\in(0,1).
\end{eqnarray}
Observe that
\begin{eqnarray*}
\dfrac{\tLambda^{X_j}}{\sum_{m=1}^d\tLambda^{X_m}}=\int\limits_{w_d(j-1)}^{w_d(j)}\dd z,\qquad \nu_d(z)=j,\quad z\in \bigl(w_d(j-1),w_d(j)\bigr],
\end{eqnarray*}
and also
\begin{eqnarray*}
	U_j(x)=\int\limits_0^{1}\id\bigl(y\leqslant U_j(x)\bigr)\dd v=\int\limits_0^{1}\id\bigl(U^{-1}_j (y)\leqslant x\bigr)\dd y,\quad x\in\R,
\end{eqnarray*}
for any $d\in\N$ and $j=1,\ldots,d$. In the latter equality we used the well known property: 
\begin{eqnarray}\label{conc_quantilefunc}
G(t)\geqslant p\quad\text{iff}\quad G^{-1}(p)\leqslant t,
\end{eqnarray} 
for any $p\in(0,1)$, $t\in\R$, and distribution function $G$ (see \cite{Vaart}, p. 304). By the above, for any $x\in\R$ and $d\in\N$ we have
\begin{eqnarray}\label{conc_Fdx}
	F_d(x)&=&\sum_{j=1}^d\Biggl(\,\,\int\limits_{w_d(j-1)}^{w_d(j)}\dd z \cdot\int\limits_0^{1}\id\bigl(U^{-1}_j (y)\leqslant x\bigr)\dd y \Biggr)\nonumber\\
	&=&\sum_{j=1}^d\int\limits_{w_d(j-1)}^{w_d(j)}\Biggl(\,\,\int\limits_0^{1}\id\bigl(U^{-1}_{\nu_d(z)} (y)\leqslant x\bigr)\dd y \Biggr)\dd z\nonumber\\
	&=&\int\limits_{0}^{1}\Biggl(\,\,\int\limits_0^{1}\id\bigl(U^{-1}_{\nu_d(z)} (y)\leqslant x\bigr)\dd y \Biggr)\dd z\nonumber\\
	&=&\iint\limits_{[0,1]^2}\id\bigl(U^{-1}_{\nu_d(z)} (y)\leqslant x\bigr)\dd y\dd z.
\end{eqnarray}

Let us consider the sequence $U^{-1}_{\nu_d(z)} (y)$, $d\in\N$. It is well known that for $n\to\infty$
\begin{eqnarray}\label{conc_quantilefunc_conv}
G_n(t)\to G(t)\,\,\text{for all}\,\, t\in\ContSet(G)\quad\text{iff}\quad G^{-1}_n(p)\to G^{-1}(p)\,\,\text{for all} \,\, p\in\ContSet(G^{-1}),
\end{eqnarray} 
where $G_n$, $n\in\N$, and $G$ are distribution functions (see \cite{Vaart} p. 305). Hence from  \eqref{th_nYd_UWconv_cond_U}  we have
\begin{eqnarray*}
	U^{-1}_j(y)-\ell_j\to U^{-1}(y),\quad j\to\infty,
\end{eqnarray*}
for all $y\in\ContSet(U^{-1})$. On account of Remark \ref{rem_sumLambda1diverg}, it is a simple matter to see that $\nu_d(z)\to\infty$ for all $z\in(0,1)$. Therefore
\begin{eqnarray}\label{conc_invGnudconv}
U^{-1}_{\nu_d(z)} (y)=U^{-1}(y)+\ell_{\nu_d(z)}+o(1),\quad d\to\infty,
\end{eqnarray}
for all $y\in\ContSet(U^{-1})$ and $z\in(0,1)$.

We now consider the sequence $\ell_{\nu_d(\cdot)}$, $d\in\N$. By the definition of $\nu_d(z)$ and the assumption of monotonicity of $(\ell_j)_{j\in\N}$, for every $z\in(0,1)$ we have
\begin{eqnarray*}
	\sum_{j=1}^d\tLambda^{X_j} \id\bigl(\signum\cdot\ell_j\leqslant \signum\cdot\ell_{\nu_d(z)}\bigr)\geqslant \sum_{j=1}^d\tLambda^{X_j}\id\bigl(j\leqslant \nu_d(z)\bigr)\geqslant z\sum_{m=1}^d\tLambda^{X_m},
\end{eqnarray*}
and also
\begin{eqnarray*}
	\sum_{j=1}^d\tLambda^{X_j}\id\bigl(\signum\cdot\ell_j< \signum\cdot\ell_{\nu_d(z)}\bigr)\leqslant \sum_{j=1}^d\tLambda^{X_j} \id\bigl(j< \nu_d(z)\bigr)< z\sum_{m=1}^d\tLambda^{X_m}.
\end{eqnarray*}
This means that $\signum\cdot\ell_{\nu_d(z)}=W^{-1}_d(z)$, $z\in(0,1)$, where
\begin{eqnarray}\label{def_Wd}
W_d(y)\defeq\dfrac{\sum_{j=1}^d\tLambda^{X_j} \id\bigl(\signum\cdot\ell_j\leqslant y\bigr)}{\sum_{m=1}^d\tLambda^{X_m}},\quad y\in\R,\quad d\in\N.
\end{eqnarray}
Since the assumption \eqref{th_nYd_UWconv_cond_W} gives the convergence $W^{-1}_d(z)-\signum\cdot a_d\to W^{-1}(z)$ for all $z\in\ContSet(W^{-1})$ as $d\to\infty$ by \eqref{conc_quantilefunc_conv}, we thus obtain
\begin{eqnarray*}
	\signum\cdot(\ell_{\nu_d(z)}-a_d)\to W^{-1}(z),\quad d\to\infty,\quad\text{for all}\quad z\in\ContSet(W^{-1}),
\end{eqnarray*}
i.e. we have
\begin{eqnarray}\label{conc_ellnudadconv}
\ell_{\nu_d(z)}-a_d\to \signum\cdot W^{-1}(z),\quad d\to\infty,\quad\text{for all}\quad z\in\ContSet(W^{-1}).
\end{eqnarray}
The relations \eqref{conc_invGnudconv} and \eqref{conc_ellnudadconv} together yield
\begin{eqnarray*}
	U^{-1}_{\nu_d(z)} (y)-a_d\to U^{-1}(y)+\signum\cdot W^{-1}(z),\quad d\to\infty,
\end{eqnarray*}
for all $y\in\ContSet(U^{-1})$ and $z\in\ContSet(W^{-1})$. From this we conclude that
\begin{eqnarray}\label{conc_idinvUinvW_conv}
\id\bigl(U^{-1}_{\nu_d(z)} (y)-a_d\leqslant x\bigr)\to \id\bigl(U^{-1}(y)+\signum\cdot W^{-1}(z)\leqslant x\bigr),\quad d\to\infty,\quad x\in\R,
\end{eqnarray}
for all $y\in\ContSet(U^{-1})$ and $z\in\ContSet(W^{-1})$ perhaps except $y$ and $z$ such that $U^{-1}(y)+\signum\cdot W^{-1}(z)=x$. Let us consider the integral
\begin{eqnarray*}
	\iint\limits_{[0,1]^2}\id\bigl(U^{-1}(y)+\signum\cdot W^{-1}(z)\leqslant x\bigr)\dd y\dd z&=&\int\limits_{[0,1]}\Biggl(\,\,\int\limits_{[0,1]}\id\bigl(U^{-1}(y)\leqslant x-\signum\cdot W^{-1}(z)\bigr)\dd y\Biggr)\dd v\\
	&=&\int\limits_{[0,1]}\Biggl(\,\,\int\limits_{[0,1]}\id\bigl(y\leqslant U(x-\signum\cdot W^{-1}(z)\bigr)\dd y\Biggr)\dd z\\
	&=&\int\limits_{[0,1]}U(x-\signum\cdot W^{-1}(z)\bigr)\dd z.
\end{eqnarray*}
Here we used the property \eqref{conc_quantilefunc}. Changing  the variables, we get
\begin{eqnarray*}
	\int\limits_{[0,1]}U(x-\signum\cdot W^{-1}(z)\bigr)\dd z=\int\limits_{-\infty}^\infty U(x-\signum\cdot v)\dd W(v)=F(x),\quad x\in\R.
\end{eqnarray*}
It is easily seen that for every $x\in \ContSet(F)$ the set $\bigl\{(y,z)\in[0,1]^2 : U^{-1}(y)+\signum\cdot W^{-1}(z)=x\bigr\}$
has zero Lebesgue measure zero. Hence the convergence \eqref{conc_idinvUinvW_conv} holds  almost everywhere on $[0,1]^2$. According to Lebesgue's dominated convergence theorem, we have
\begin{eqnarray*}
	\iint\limits_{[0,1]^2}\id\bigl(U^{-1}_{\nu_d(z)} (y)-a_d\leqslant x\bigr)\dd y\dd z\to\iint\limits_{[0,1]^2}\id\bigl(U^{-1}(y)+\signum\cdot W^{-1}(z)\leqslant x\bigr)\dd y\dd z,\quad d\to\infty,
\end{eqnarray*}
i.e. we prove \eqref{conc_FdtoFad}.

Next, we consider the integral from \eqref{conc_nYdesimint}
\begin{eqnarray*}
	\int\limits_{0}^{1-(\e/\e_d)^2} \exp\bigl\{F^{-1}_d(y)\bigr\}\dd y=e^{a_d}\int\limits_{0}^{1} B_d(y)\dd y,
\end{eqnarray*}
where 
\begin{eqnarray*}
	B_d(y)\defeq \exp\bigl\{F^{-1}_d(y)-a_d\bigr\}\id\bigl(y\leqslant 1-(\e/\e_d)^2\bigr),\quad y\in(0,1), \quad d\in\N.
\end{eqnarray*}
By \eqref{conc_quantilefunc_conv} the convergence \eqref{conc_FdtoFad} yield 
\begin{eqnarray}\label{conc_invFdconv}
F^{-1}_d(y)-a_d\to F^{-1}(y),\quad d\to\infty,\quad\text{for all}\quad y\in\ContSet(F^{-1}).
\end{eqnarray}
On account of monotonicity of $F^{-1}$, the set $(0,1)\setminus \ContSet(F^{-1})$ is countable. Hence we have almost everywhere convergence of $F^{-1}_d(y)-a_d$ to $F^{-1}(y)$ on $(0,1)$ as $d\to\infty$. This yields 
\begin{eqnarray*}
	B_d(y)\to\exp\bigl\{F^{-1}(y)\bigr\}\id\bigl(y\leqslant 1-(\e/\e_0)^2\bigr)\quad\text{for almost all}\quad y\in(0,1).
\end{eqnarray*}
The functions $B_d$, $d\in\N$, are uniformly bounded on $(0,1)$. Indeed, let us fix $y_\e\in \ContSet(F^{-1})$ such that   $ 1-(\e/\e_0)^2\leqslant y_\e<1$. Since $F^{-1}_d$, $d\in\N$, are non-decreasing, we have 
\begin{eqnarray*}
	B_d(y)\leqslant \exp\bigl\{F^{-1}_d(y_\e) -a_d\bigr\}\quad\text{for all}\quad y\in(0,1), \quad d\in\N.
\end{eqnarray*}
By the relation \eqref{conc_invFdconv}, the sequence $F^{-1}_d(y_\e) -a_d$ converges to $F^{-1}(y_\e)$, and hence bounded. This gives uniform boundedness of $B_d$, $d\in\N$.
According to Lebesgue's dominated convergence theorem, we have
\begin{eqnarray*}
	\int\limits_{0}^{1} B_d(y)\dd y\to \int\limits_{0}^{1}\Bigl(\exp\bigl\{F^{-1}(y)\bigr\}\id\bigl(y\leqslant 1-(\e/\e_0)^2\bigr)\Bigr)\dd y,\quad d\to\infty.
\end{eqnarray*}
It exactly means that
\begin{eqnarray*}
	\int\limits_{0}^{1-(\e/\e_d)^2} \exp\bigl\{F^{-1}_d(y)\bigr\}\dd y\sim e^{a_d}\!\!\!\!\! \int\limits_{0}^{1-(\e/\e_0)^2} \exp\bigl\{F^{-1}(y)\bigr\}\dd y,\quad d\to\infty.
\end{eqnarray*}
From this and \eqref{conc_nYdesimint} we come to \eqref{th_nYd_UWconv_conc}. \quad$\Box$\\

In some cases there is not any centering sequence $(a_d)_{d\in\N}$ to satisfy \eqref{th_nYd_UWconv_cond_W}. Often this problem is solved by introducing admissible norming sequence $(b_d)_{d\in\N}$, which tends to infinity.  For these cases we have the following result.

\begin{Theorem}\label{th_nYd_UWconv_log}
	Let $X_j$, $j\in\N$, satisfy Assumptions $\ref{assum_Regular}$ and $\ref{assum_tendstoinfty}$. 
	Let $\e_0\in(0,1]$ be the number from $\eqref{def_e0}$. Suppose that
	\begin{eqnarray*}
	\dfrac{1  }{\tLambda^{X_j}}\sum_{k=1}^\infty\tlambda^{X_j}_k\,\id\bigl(-\ln\tlambda^{X_j}_k-\ell_j\leqslant x\bigr)\to U(x),\quad j\to\infty,\quad\text{for all}\quad x\in\ContSet(U),
	\end{eqnarray*}
	with a distribution function $U$ and a monotonic sequence $(\ell_j)_{j\in\N}$.
	Let $\signum=1$ if $(\ell_j)_{j\in\N}$ is non-decreasing, and  $\signum=-1$ if $(\ell_j)_{j\in\N}$ is non-increasing. Next, suppose that
	\begin{eqnarray}\label{th_nYd_UWconv_log_cond_W}
	\dfrac{\sum_{j=1}^d\tLambda^{X_j} \id\bigl(\signum\cdot(\ell_j-a_d)\leqslant x b_d\bigr)}{\sum_{m=1}^d\tLambda^{X_m}}\to W(x),\quad d\to\infty,\quad\text{for all}\quad x\in\ContSet(W),
	\end{eqnarray}
	with a distribution function $W$, a sequence $(a_d)_{d\in\N}$ and a positive sequence $(b_d)_{d\in\N}$ such that $b_d\to\infty$, $d\to\infty$.
	Then the following asymptotics holds
	\begin{eqnarray}\label{th_nYd_UWconv_log_conc}
		\ln n^{Y_d}(\e)=\ln \Bigl(\sum_{j=1}^d \tLambda^{X_j}\Bigr) +a_d +F^{-1}\bigl(1-(\e/\e_0)^2\bigr) b_d+o(b_d),\quad d\to\infty,
	\end{eqnarray}
	for every $\e\in(0,\e_0)$ such that $1-(\e/\e_0)^2\in\ContSet(F^{-1})$. Here $F(x)=W(x)$, $x\in\R$, if $\signum=1$, and $F(x)=1-W(x-0)$, $x\in\R$, if $\signum=-1$. 
\end{Theorem}
\textbf{Proof.}\quad Let $\e_d$ and $F_d$, $d\in\N$, be defined by \eqref{def_ed} and \eqref{def_Fd} respectively. By Assumption \ref{assum_Regular}, we have $\e_d\to\e_0$ as $d\to\infty$. Let us choose $\e\in(0,\e_0)$ such that $p_\e\defeq 1-(\e/\e_0)^2$ belongs to $\ContSet(F^{-1})$. Thus Theorem~\ref{th_nYd_repr} yields the representation \eqref{th_nYd_repr_conc} of $n^{Y_d}(\e)$ for all sufficiently large $d\in\N$.

Let us consider the sequence of distribution functions $F_d$, $d\in\N$. We first show that
\begin{eqnarray}\label{conc_FdtoFadbd}
F_d(xb_d+a_d)\to F(x),\quad d\to\infty,\quad \text{for all}\quad x\in\ContSet(F).
\end{eqnarray}
We define the functions $U_j$, $\nu_d$ and $W_d$ by formulas \eqref{def_Uj}, \eqref{def_wd_nud} and \eqref{def_Wd} respectively. According to \eqref{conc_quantilefunc_conv}, the assumption \eqref{th_nYd_UWconv_log_cond_W} gives the convergence
\begin{eqnarray*}
	\dfrac{W^{-1}_d(z)-\signum\cdot a_d}{b_d}\to W^{-1}(z),\quad d\to\infty,\quad\text{for all}\quad z\in\ContSet(W^{-1}).
\end{eqnarray*}
In the proof of Theorem \ref{th_nYd_UWconv} it was shown that $\signum\cdot\ell_{\nu_d(z)}=W^{-1}_d(z)$, $z\in(0,1)$. Hence
\begin{eqnarray*}
\dfrac{\ell_{\nu_d(z)}-a_d}{b_d}\to\signum\cdot W^{-1}(z),\quad d\to\infty,\quad\text{for all}\quad z\in\ContSet(W^{-1}).
\end{eqnarray*}
On account of the proved asymptotic relation \eqref{conc_invGnudconv}, we have
\begin{eqnarray*}
	\dfrac{	U^{-1}_{\nu_d(z)} (y)-a_d}{b_d}=\signum\cdot W^{-1}(z)+ U^{-1}(y)/b_d+o(1/b_d),\quad d\to\infty.
\end{eqnarray*}
for all $y\in\ContSet(U^{-1})$ and $z\in\ContSet(W^{-1})$. Since $b_d\to\infty$, $d\to\infty$, by the assumption, we obtain for these $y$ and $z$
\begin{eqnarray}\label{conc_Uinvnudzyadbd_conv}
\dfrac{	U^{-1}_{\nu_d(z)} (y)-a_d}{b_d}\to\signum\cdot W^{-1}(z),\quad d\to\infty.
\end{eqnarray}
From the equality \eqref{conc_Fdx} we get 
\begin{eqnarray*}
	F_d(x b_d+a_d)=	\iint\limits_{[0,1]^2}\id\bigl(U^{-1}_{\nu_d(z)} (y)-a_d\leqslant x b_d\bigr)\dd y\dd z,\quad x\in\R,\quad d\in\N.
\end{eqnarray*}
Now we consider the integral
\begin{eqnarray*}
	\iint\limits_{[0,1]^2}\id\bigl(\signum\cdot W^{-1}(z)\leqslant x\bigr)\dd y\dd z=\int\limits_{[0,1]}\id\bigl(\signum\cdot W^{-1}(z)\leqslant x\bigr)\dd z
	=\int\limits_{-\infty}^\infty\id\bigl(\signum\cdot v\leqslant x\bigr)\dd W(v).
\end{eqnarray*}
It is easy to check that the last integral is equal to $F(x)$ for each $x\in\R$. Let us fix $x\in\ContSet(F)$. Then the set $\bigl\{(y,z)\in[0,1]^2:\signum\cdot W^{-1}(z)=x\bigr\}$ has zero Lebesgue measure. This fact and \eqref{conc_Uinvnudzyadbd_conv} together imply the convergence
\begin{eqnarray*}
	\id\bigl(U^{-1}_{\nu_d(z)} (y)-a_d\leqslant x b_d\bigr)\to \id\bigl(\signum\cdot W^{-1}(z)\leqslant x\bigr),\quad d\to\infty,
\end{eqnarray*}
which holds for almost all $(y,z)\in[0,1]^2$. Next, using Lebesgue's dominated convergence theorem  we conclude that
\begin{eqnarray*}
	\iint\limits_{[0,1]^2}\id\bigl(U^{-1}_{\nu_d(z)} (y)-a_d\leqslant x b_d\bigr)\dd y\dd z\to\iint\limits_{[0,1]^2}\id\bigl(\signum\cdot W^{-1}(z)\leqslant x\bigr)\dd y\dd z,\quad d\to\infty.
\end{eqnarray*}
Thus we proved the convergence \eqref{conc_FdtoFadbd}.

Next, we return to the formula \eqref{th_nYd_repr_conc}. On account of Assumption \ref{assum_tendstoinfty} and Proposition \ref{pr_nde_toinfty}, we have $n^{Y_d}(\e)\to\infty$ as $d\to\infty$. Fix arbitrarily small $\delta>0$. On account of $b_d\to\infty$, $d\to\infty$, we conclude that
\begin{eqnarray}\label{conc_ndeineqright}
n^{Y_d}(\e)\leqslant\sum_{j=1}^d \tLambda^{X_j}\cdot\!\!\! \int\limits_{0}^{1-(\e/\e_d)^2} \exp\bigl\{F^{-1}_d(y)\bigr\}\dd y\cdot e^{\delta b_d},
\end{eqnarray}
for all  large enough  $d\in\N$.
Since every $F^{-1}_d$ is non-decreasing, we have
\begin{eqnarray}\label{conc_intexpFdinveedright}
\int\limits_{0}^{1-(\e/\e_d)^2} \exp\bigl\{F^{-1}_d(y)\bigr\}\dd x\leqslant\exp\bigl\{F^{-1}_d\bigl(1-(\e/\e_d)^2\bigr)\bigr\}.
\end{eqnarray}
The function $F^{-1}$ is non-decreasing, and consequently the set $\ContSet(F^{-1})$ is dense in the interval $(0,1)$. Hence we can choose $\tau_1>0$ such that $p_\e+\tau_1\in \ContSet(F^{-1})$ and $F^{-1}(p_\e+\tau_1)-F^{-1}(p_\e)\leqslant \delta$. Since, by \eqref{conc_quantilefunc_conv}, the convergence \eqref{conc_FdtoFadbd} yields
\begin{eqnarray}\label{conc_Finvdadbdconv}
\dfrac{F^{-1}_d(p)-a_d}{b_d}\to F^{-1}(p),\quad d\to\infty,\quad\text{for all}\quad p\in\ContSet(F^{-1}),
\end{eqnarray}
in particular, we have
\begin{eqnarray*}
	\dfrac{F^{-1}_d(p_\e+\tau_1)-a_d}{b_d}\to F^{-1}(p_\e+\tau_1),\quad d\to\infty.
\end{eqnarray*}
Hence for all large enough $d\in\N$ we obtain
\begin{eqnarray*}
	F^{-1}_d\bigl(1-(\e/\e_d)^2\bigr)&\leqslant& F^{-1}_d(p_\e+\tau_1)\\
	&\leqslant& a_d+ F^{-1}(p_\e+\tau_1)b_d+\delta b_d\\
	&\leqslant& a_d+ F^{-1}(p_\e)b_d+2\delta b_d.
\end{eqnarray*}
Combining this inequality with \eqref{conc_ndeineqright} and \eqref{conc_intexpFdinveedright} we obtain the following inequality
\begin{eqnarray*}
	\ln n^{Y_d}(\e)\leqslant \ln\Bigl(\sum_{j=1}^d \tLambda^{X_j}\Bigr)+a_d+ F^{-1}(p_\e)b_d+3\delta b_d,
\end{eqnarray*}
which holds for all large enough $d\in\N$.

We now get a similar lower bound for $\ln n^{Y_d}(\e)$. From the formula \eqref{th_nYd_repr_conc} we directly have
\begin{eqnarray}\label{conc_ndeineqleft}
n^{Y_d}(\e)\geqslant\sum_{j=1}^d \tLambda^{X_j}\cdot\!\!\! \int\limits_{0}^{1-(\e/\e_d)^2} \exp\bigl\{F^{-1}_d(y)\bigr\}\dd y,
\end{eqnarray}
for all sufficiently large $d\in\N$. Let us choose $\tau_2>0$ such that $p_\e-\tau_2\in \ContSet(F^{-1})$ and $F^{-1}(p_\e)-F^{-1}(p_\e-\tau_2)\leqslant \delta$. Since $\e_d\to\e_0$, $d\to\infty$,  we have the following inequality for all large enough $d\in\N$  
\begin{eqnarray*}
	1-(\e/\e_d)^2-(p_\e-\tau_2)\geqslant \tau_2/2>0.
\end{eqnarray*}
For these $d$ the integral from  \eqref{conc_ndeineqleft} we estimate in the following way
\begin{eqnarray}\label{conc_intexpFdinveedleft}
\int\limits_{0}^{1-(\e/\e_d)^2} \exp\bigl\{F^{-1}_d(y)\bigr\}\dd x&\geqslant& \int\limits_{p_\e-\tau_2}^{1-(\e/\e_d)^2} \exp\bigl\{F^{-1}_d(y)\bigr\}\dd x\nonumber\\
&\geqslant &\exp\bigl\{F^{-1}_d(p_\e-\tau_2)\bigr\}\bigl(1-(\e/\e_d)^2-(p_\e-\tau_2)\bigr)\nonumber\\
&\geqslant &\exp\bigl\{F^{-1}_d(p_\e-\tau_2)\bigr\}\tau_2/2.
\end{eqnarray}
From \eqref{conc_Finvdadbdconv} we have
\begin{eqnarray*}
	\dfrac{F^{-1}_d(p_\e-\tau_2)-a_d}{b_d}\to F^{-1}(p_\e-\tau_2),\quad d\to\infty.
\end{eqnarray*}
Consequently, for all large enough $d\in\N$ we have
\begin{eqnarray*}
	F^{-1}_d(p_\e-\tau_2)&\geqslant& a_d+ F^{-1}(p_\e-\tau_2)b_d-\delta b_d\\
	&\geqslant& a_d+ F^{-1}(p_\e)b_d-2\delta b_d.
\end{eqnarray*}
Combining the latter inequality with \eqref{conc_ndeineqleft} and \eqref{conc_intexpFdinveedleft} gives
\begin{eqnarray*}
	\ln n^{Y_d}(\e)\geqslant \ln\Bigl(\sum_{j=1}^d \tLambda^{X_j}\Bigr)+a_d+ F^{-1}(p_\e)b_d-2\delta b_d+\ln(\tau_2/2)
\end{eqnarray*}
for all large enough $d\in\N$. On account of $b_d\to\infty$, $d\to\infty$, we obtain
the inequality
\begin{eqnarray*}
	\ln n^{Y_d}(\e)\geqslant \ln\Bigl(\sum_{j=1}^d \tLambda^{X_j}\Bigr)+a_d+ F^{-1}(p_\e)b_d-3\delta b_d,
\end{eqnarray*}
which holds for all sufficiently large $d\in\N$.

The obtained upper and lower estimates of $\ln n^{Y_d}(\e)$ yield the required asymptotics \eqref{th_nYd_UWconv_log_conc}.\quad $\Box$

\section{Applications to Korobov kernels}

Let $B_{\alpha,\beta,\sigma}(t)$, $t\in[0,1]$, be a zero-mean random process with the following covariance function
\begin{eqnarray*}
	\CovFunc_{\alpha,\beta,\sigma}(t,s)\colonequals \alpha+2\beta\sum_{k=1}^{\infty} k^{-\sigma}\cos(2\pi k(t-s)), \quad t,s\in[0,1],
\end{eqnarray*} 
which is called \textit{Korobov kernel}. Here $\alpha>0$, $\beta>0$ and $\sigma>1$. 

Let us recall eigenpairs of the covariance operator $K^{B_{\alpha,\beta,\sigma}}$ (see \cite{NovWoz1}, Appendix A). The identical $1$ is an eigenvector of $K^{B_{\alpha,\beta,\sigma}}$ with the eigenvalue $\tlambda^{B_{\alpha,\beta,\sigma}}_0=\alpha$. The other eigenpairs of $K^{B_{\alpha,\beta,\sigma}}$ have the following form:
\begin{eqnarray*} 
\tlambda^{B_{\alpha,\beta,\sigma}}_{2k-1}=\tlambda^{B_{\alpha,\beta,\sigma}}_{2k}=\dfrac{\beta}{k^{\sigma}},\quad \tpsi^{B_{\alpha,\beta,\sigma}}_{2k-1}(t)=e^{-i 2\pi k t}, \quad \tpsi^{B_{\alpha,\beta,\sigma}}_{2k}(t)=e^{i 2\pi k t},\quad k\in\N,\quad t\in[0,1].
\end{eqnarray*}
Note that the trace of $K^{B_{\alpha,\beta,\sigma}}$ is
\begin{eqnarray*}
	\Lambda^{B_{\alpha,\beta,\sigma}}=\alpha+2\beta \zeta(\sigma),
\end{eqnarray*}
where $\zeta(p)=\sum_{k=1}^\infty k^{-p}$, $p>1$, is the Riemann zeta-function.

Suppose that we have a sequence of processes $B_j(t)$, $t\in[0,1]$,  with covariance functions $\CovFunc_{\alpha_j,\beta_j,\sigma_j}$, $j\in\N$, respectively. Let $\KorElem_d(t)$, $t\in[0,1]^d$, $d\in\N$, be the sequence of zero-mean random fields with the following covariance functions
\begin{eqnarray*}
	\CovFunc^{\KorElem_d}(t,s)=\sum_{j=1}^{d}\CovFunc_{\alpha_j,\beta_j,\sigma_j}(t_j,s_j),
\end{eqnarray*}
where $t=(t_1,\ldots, t_d)$ and $s=(s_1,\ldots, s_d)$ are from $[0,1]^d$, $d\in\N$.
We consider every field $\KorElem_d(t)$, $t\in [0,1]^d$, as a random element of the space $L_2([0,1]^d)$.  We investigate asymptotic behaviour of the approximation complexity $n^{\KorElem_d}(\e)$ for fixed $\e\in(0,1)$ and $d\to\infty$. Our basic assumption holds for the sequence of marginal processes $B_j(t)$, $t\in[0,1]$, $j\in\N$.

In order to evidently illustrate the application of the general results from previous sections and to avoid routine unwieldy calculations, we will solve our approximation problem under the following assumptions on the parameters:
\begin{eqnarray}\label{assum_param}
\beta_j\sim c j^{-\tau},\qquad \alpha_j/\beta_j\to r,\qquad\sigma_j\to\infty,\qquad  j\to\infty,
\end{eqnarray}
where $c>0$, $\tau\in\R$, and  $0\leqslant r\leqslant \infty$.

\begin{Proposition}
Let $B_j$, $j\in\N$, satisfy  $\eqref{assum_param}$ with some $c>0$ and either $\tau>1$, $0\leqslant r\leqslant \infty$, or $\tau\leqslant 1$, $r=\infty$. Then 
\begin{eqnarray*}
\sup_{d\in\N} n^{\KorElem_d}(\e)<\infty,\quad \text{for every}\quad \e\in(0,1).
\end{eqnarray*}
\end{Proposition}
\textbf{Proof.}\quad Suppose that $\tau>1$ and $0\leqslant r\leqslant \infty$. Let us consider the series
\begin{eqnarray*}
\sum_{j=1}^\infty \tLambda^{B_j}=\sum_{j=1}^\infty 2\beta_j \zeta(\sigma_j).
\end{eqnarray*}
Since, by the assumptions, $\sum_{j=1}^\infty \beta_j<\infty$ and $\zeta(\sigma_j)\to1$ as $j\to \infty$, we have the convergence of $\sum_{j=1}^\infty \tLambda^{B_j}$. Using Proposition \ref{pr_nYdbounded} we get the required assertion.

Suppose that $\tau\leqslant1$ and $r=\infty$. Therefore
\begin{eqnarray*}
\sum_{j=1}^d 2\beta_j \zeta(\sigma_j)=o\Bigl(\sum_{j=1}^d \alpha_j\Bigr),\quad d\to\infty.
\end{eqnarray*}
Consequently, the sequence
\begin{eqnarray*}
\dfrac{\sum_{j=1}^d\tlambda^{B_j}_0}{\sum_{j=1}^d\Lambda^{B_j}}=\dfrac{\sum_{j=1}^d\alpha_j}{\sum_{j=1}^d\alpha_j+ \sum_{j=1}^{d} 2\beta_j \zeta(\sigma_j)},\quad d\in\N,
\end{eqnarray*}
tends to $1$ as $d\to\infty$. By the remarks before Assumption \ref{assum_Regular}, we have $n^{\KorElem_d}(\e)=1$ for all large enough $d\in\N$. Obviously, this implies the  boundedness of $n^{\KorElem_d}(\e)$ on $d\in\N$ for every $\e\in(0,1)$. \quad $\Box$ 

Now we focus only on the case $c>0$, $\tau\leqslant 1$, $0\leqslant r<\infty$. Here  $B_j$, $j\in\N$, satisfy Assumption $\ref{assum_Regular}$. Indeed, since
$\sum_{j=1}^d \beta_j\to\infty$ and $\zeta(\sigma_d)\to1$, $d\to\infty$, we have
\begin{eqnarray*}
	\dfrac{\sum_{j=1}^{d} \beta_j \zeta(\sigma_j)}{\sum_{j=1}^d\beta_j}\to 1,\quad d\to\infty,
\end{eqnarray*}
and also
\begin{eqnarray*}
	\dfrac{\sum_{j=1}^{d} \alpha_j }{\sum_{j=1}^d\beta_j}=\dfrac{\sum_{j=1}^{d} (\alpha_j/\beta_j)\beta_j }{\sum_{j=1}^d\beta_j}\to r,\quad d\to\infty,
\end{eqnarray*}
where we used the well known theorems of summability of numerical series (see \cite{Zygmund} p. 74). So we obtain
\begin{eqnarray*}
	\dfrac{\sum_{j=1}^d\tlambda^{B_j}_0}{\sum_{j=1}^d\Lambda^{B_j}}=\dfrac{\sum_{j=1}^d\alpha_j}{\sum_{j=1}^d\alpha_j+ \sum_{j=1}^{d} 2\beta_j \zeta(\sigma_j)}\to\dfrac{r}{r+2},\quad d\to\infty.
\end{eqnarray*}
Thus Assumption $\ref{assum_Regular}$ holds and  $\e_0=(1+r/2)^{-1/2}\leqslant 1$ by \eqref{def_e0}. According to the remarks after Assumption \ref{assum_Regular}, it makes sense to consider $n^{\KorElem_d}(\e)$ only for $\e\in(0,\e_0)$ if $\e_0<1$. 

For $\tau<1$ and $\tau=1$ the approximation complexity grows rather differently. So we consider these cases separately.

\begin{Theorem}\label{th_nKorrElemde}
Let $B_j$, $j\in\N$, satisfy  $\eqref{assum_param}$ with $c>0$, $\tau<1$, and  $0\leqslant r<\infty$. Then
\begin{eqnarray*}
	n^{\KorElem_d}(\e)\sim 2Q(\e)\cdot d ,\quad d\to\infty,\quad\text{for every}\quad \e\in(0,\e_0),
\end{eqnarray*}
where
\begin{eqnarray*}
Q(\e)=
\begin{cases}
	\bigl(1-(\e/\e_0)^2\bigr)^{\tfrac{1}{1-\tau}},&\text{if $\tau\in[0,1)$},\\
	1-(\e/\e_0)^{\tfrac{2}{1-\tau}},& \text{if   $\tau<0$},
\end{cases}
\qquad \e\in(0,\e_0).
\end{eqnarray*}
\end{Theorem}
\textbf{Proof.}\quad We first check that $B_j$, $j\in\N$, satisfy Assumption $\ref{assum_tendstoinfty}$. Let us consider the sequence
\begin{eqnarray*}
	\dfrac{\max\limits_{j=1,\ldots,d}\tlambda^{B_j}_1}{\sum_{j=1}^d \tLambda^{B_j}}=\dfrac{\max\limits_{j=1,\ldots,d} \beta_j}{2\sum_{j=1}^d \beta_j \zeta(\sigma_j)},\quad d\in\N.
\end{eqnarray*}
On account of the assumption for $(\beta_j)_{j\in\N}$, it is easy to check that
for some constant $C>0$
\begin{eqnarray*}
\max\limits_{j=1,\ldots,d} \beta_j\leqslant C \max\{1,d^{-\tau}\}\quad \text{for all}\quad d\in\N.
\end{eqnarray*}
Since $\sum_{j=1}^d \beta_j\to\infty$ and $\zeta(\sigma_d)\to1$, $d\to\infty$, the following equivalences hold
\begin{eqnarray*}
\sum_{j=1}^d \beta_j \zeta(\sigma_j)\sim \sum_{j=1}^d \beta_j \sim \sum_{j=1}^d c j^{-\tau},\quad d\to\infty.
\end{eqnarray*}
Hence we have
\begin{eqnarray}\label{conc_betajzetasjtauless1}
\sum_{j=1}^d \beta_j \zeta(\sigma_j)\sim  \dfrac{c\, d^{1-\tau}}{1-\tau},\quad d\to\infty.
\end{eqnarray}
We see that
\begin{eqnarray*}
\dfrac{\max\limits_{j=1,\ldots,d}\tlambda^{B_j}_1}{\sum_{j=1}^d \tLambda^{B_j}}=O(1/d),\quad d\to\infty.
\end{eqnarray*}
Thus this sequence tends to zero, and the condition from Assumption $\ref{assum_tendstoinfty}$ is true.

Let us define $\ell_j\defeq -\ln c+\tau \ln j$, $j\in\N$, and consider the sums from \eqref{th_nYd_UWconv_cond_U} for $X_j=B_j$, $j\in\N$, $x\in\R$:
\begin{eqnarray*}
\dfrac{1  }{\tLambda^{B_j}}\sum_{k=1}^\infty\tlambda^{B_j}_k\,\id\bigl(-\ln\tlambda^{B_j}_k-\ell_j\leqslant x\bigr)
&=&\dfrac{1  }{\tLambda^{B_j}}\sum_{k=1}^\infty(\tlambda^{B_j}_{2k-1}+\tlambda^{B_j}_{2k})\,\id\bigl(-\ln\tlambda^{B_j}_{2k}-\ell_j\leqslant x\bigr)\\
&=&\dfrac{1  }{2\beta_j\zeta(\sigma_j)}\sum_{k=1}^\infty\dfrac{2\beta_j}{k^{\sigma_j}}\,\id\bigl(\sigma_j\ln k-\ln \beta_j-\ell_j\leqslant x\bigr)\\
&=&\dfrac{1  }{\zeta(\sigma_j)}\sum_{k=1}^\infty\dfrac{1}{k^{\sigma_j}}\,\id\bigl(\sigma_j\ln k+\ln(cj^{-\tau}/\beta_j)\leqslant x\bigr).
\end{eqnarray*}
Since $\sigma_j\to\infty$ and $\ln(cj^{-\tau}/\beta_j) =o(1)$ as $j\to\infty$, we obtain the convergence \eqref{th_nYd_UWconv_cond_U} with $U(x)=\id(x\geqslant 0)$, $x\in\R$.

Next, we consider the sequence from left-hand side of \eqref{th_nYd_UWconv_cond_W} with $X_j=B_j$, $j\in\N$, and  $a_d\defeq \ell_d $, $d\in\N$. According to monotonicity of $(\ell_j)_{j\in\N}$, we set $\signum\defeq \sign(\tau)$ if $\tau\ne 0$, and $\signum\defeq 1$ if $\tau=0$. For any $x\in\R$ and $d\in\N$ we have
\begin{eqnarray*}
\dfrac{\sum_{j=1}^d\tLambda^{B_j} \id\bigl(\signum\cdot(\ell_j-a_d)\leqslant x\bigr)}{\sum_{j=1}^d\tLambda^{B_j}}
&=& \dfrac{\sum_{j=1}^d 2\beta_j \zeta(\sigma_j) \id\bigl(\signum\cdot(\ell_j-\ell_d)\leqslant x\bigr)}{\sum_{j=1}^d2\beta_j \zeta(\sigma_j)}\\
&=& \dfrac{\sum_{j=1}^d \beta_j \zeta(\sigma_j) \id\bigl(\signum\cdot\ln(j/d)^\tau\leqslant x\bigr)}{\sum_{j=1}^d \beta_j \zeta(\sigma_j)}\\
&=& \dfrac{\sum_{j=1}^d \beta_j \zeta(\sigma_j) \id\bigl(|\tau|\cdot\ln(j/d)\leqslant x\bigr)}{\sum_{j=1}^d \beta_j \zeta(\sigma_j)}.
\end{eqnarray*}
If $\tau=0$, then this fraction exactly equals $\id(x\geqslant 0)$ for all $x\in\R$. Hence \eqref{th_nYd_UWconv_cond_W} holds with $W(x)=\id(x\geqslant 0)$, $x\in\R$. If $\tau\ne 0$, then we can write
 \begin{eqnarray*}
 \dfrac{\sum_{j=1}^d\tLambda^{B_j} \id\bigl(\signum\cdot(\ell_j-a_d)\leqslant x\bigr)}{\sum_{j=1}^d\tLambda^{B_j}}= \dfrac{\sum_{j=1}^d \beta_j \zeta(\sigma_j) \id\bigl(j\leqslant d\cdot e^{x/|\tau|}\bigr)}{\sum_{j=1}^d \beta_j \zeta(\sigma_j)},\quad x\in\R,\quad d\in\N.
 \end{eqnarray*}
For $x\geqslant 0$ the latter fraction is equal to identically $1$. For $x<0$ and all large enough $d$ we have
\begin{eqnarray*}
\sum_{j=1}^d \beta_j \zeta(\sigma_j) \id\bigl(j\leqslant d\cdot e^{x/|\tau|}\bigr)=\sum_{j=1}^{j_{d,x}} \beta_j \zeta(\sigma_j),
\end{eqnarray*}
where $j_{d,x}\defeq \max\bigl\{j\in\N:j\leqslant d\cdot e^{x/|\tau|}\bigr\}$. It is  a simple matter to check that $j_{d,x}\sim d\cdot e^{x/|\tau|}$, $d\to\infty$. 
Using this and the relation \eqref{conc_betajzetasjtauless1} we obtain 
\begin{eqnarray*}
\dfrac{\sum_{j=1}^{j_{d,x}} \beta_j \zeta(\sigma_j)}{\sum_{j=1}^d \beta_j \zeta(\sigma_j)}\sim\biggl(\dfrac{j_{d,x}}{d}\biggr)^{1-\tau}\to \exp\bigl\{\tfrac{1-\tau} {|\tau|}\,x\bigr\},\quad d\to\infty.
\end{eqnarray*}
Thus the convergence \eqref{th_nYd_UWconv_cond_W} holds with 
\begin{eqnarray*}
W(x)=\exp\bigl\{\tfrac{1-\tau} {|\tau|}\,x\bigr\}\id(x<0)+\id(x\geqslant 0),\quad x\in\R.
\end{eqnarray*}

All conditions of Theorem \ref{th_nYd_UWconv} hold. According to this theorem,  for every $\e\in(0,\e_0)$ we have 
\begin{eqnarray}\label{conc_nBde}
n^{\KorElem_d}(\e)\sim e^{\ell_d}\sum_{j=1}^d 2\beta_j\zeta(\sigma_j)\!\!\!\!\! \int\limits_{0}^{1-(\e/\e_0)^2} \exp\bigl\{F^{-1}(y)\bigr\}\dd y,\qquad d\to\infty,
\end{eqnarray}
where $F(x)=\int_{-\infty}^\infty U(x-\signum\cdot v)\dd W(v)$, $x\in\R$. Since $e^{\ell_d}=d^\tau/c$, $d\in\N$, using \eqref{conc_betajzetasjtauless1} we conclude
\begin{eqnarray*}
	e^{\ell_d}\sum_{j=1}^d 2\beta_j\zeta(\sigma_j)\sim  \dfrac{2\, d}{1-\tau},\quad d\to\infty.
\end{eqnarray*}
Since $U(x)=\id(x\geqslant 0)$, $x\in\R$, we have  $F(x)=\int_{\signum\cdot v\leqslant x} \dd W(v)$, $x\in\R$. Hence $F(x)=W(x)$, $x\in\R$, if $\signum=1$, and $F(x)=1-W(-x-0)$, $x\in\R$, if $\signum=-1$. 

Separating the cases $\tau=0$, $\tau\in(0,1)$ and $\tau<0$,  we now obtain explicit formulas
for the right-hand side of \eqref{conc_nBde}. In the case $\tau=0$ we have $\signum=1$ and $F(x)=W(x)=\id(x\geqslant 0)$, $x\in\R$. Hence $F^{-1}(y)=1$, $y\in(0,1)$, and
\begin{eqnarray*}
\int\limits_{0}^{1-(\e/\e_0)^2} \exp\bigl\{F^{-1}(y)\bigr\}\dd y=\int\limits_{0}^{1-(\e/\e_0)^2}\dd y=1-(\e/\e_0)^2.
\end{eqnarray*}
Thus we obtain
\begin{eqnarray*}
n^{\KorElem_d}(\e)\sim 2\,\bigl(1-(\e/\e_0)^2\bigr)\cdot d ,\qquad d\to\infty,\quad \e\in(0,\e_0).
\end{eqnarray*}
In the case  $\tau\in(0,1)$ we have $\signum=1$ and
\begin{eqnarray*}
F(x)=W(x)=\exp\bigl\{\tfrac{1-\tau} {\tau}\,x\bigr\}\id(x<0)+\id(x\geqslant 0),\quad x\in\R.
\end{eqnarray*}
Consequently, $F^{-1}(y)=\tfrac{\tau}{1-\tau}\ln y$, $y\in(0,1)$, and
\begin{eqnarray*}
	\int\limits_{0}^{1-(\e/\e_0)^2} \exp\bigl\{F^{-1}(y)\bigr\}\dd y=\int\limits_{0}^{1-(\e/\e_0)^2} y^{\tfrac{\tau}{1-\tau}}\dd y=(1-\tau)\bigl(1-(\e/\e_0)^2\bigr)^{\tfrac{1}{1-\tau}}.
\end{eqnarray*}
Hence we obtain
\begin{eqnarray*}
	n^{\KorElem_d}(\e)\sim 2\,\bigl(1-(\e/\e_0)^2\bigr)^{\tfrac{1}{1-\tau}}\cdot d ,\qquad d\to\infty,\quad \e\in(0,\e_0).
\end{eqnarray*}
In the case $\tau<0$ we have $\signum=-1$, $W(x)=\exp\bigl\{-\tfrac{1-\tau} {\tau}\,x\bigr\}\id(x<0)+\id(x\geqslant 0)$, $x\in\R$, and
\begin{eqnarray*}
	F(x)=1-W(-x-0)=\bigl(1-\exp\bigl\{\tfrac{1-\tau} {\tau}\,x\bigr\}\bigr)\id(x\geqslant 0),\quad x\in\R.
\end{eqnarray*}
Therefore $F^{-1}(y)=\tfrac{\tau}{1-\tau}\ln (1-y)$, $y\in(0,1)$, and
\begin{eqnarray*}
	\int\limits_{0}^{1-(\e/\e_0)^2} \exp\bigl\{F^{-1}(y)\bigr\}\dd y=\int\limits_{0}^{1-(\e/\e_0)^2} (1-y)^{\tfrac{\tau}{1-\tau}}\dd y=(1-\tau)\biggl(1-(\e/\e_0)^{\tfrac{2}{1-\tau}}\biggr).
\end{eqnarray*}
Hence we  have 
\begin{eqnarray*}
	n^{\KorElem_d}(\e)\sim 2\,\biggl(1-(\e/\e_0)^{\tfrac{2}{1-\tau}}\biggr)\cdot d ,\qquad d\to\infty, \quad \e\in(0,\e_0).
\end{eqnarray*}
Thus we obtained the required asymptotics for $n^{\KorElem_d}(\e)$. \quad $\Box$\\

Under the assumption $\beta_j\sim c/j$, $j\to\infty$, without any supplements, we can not obtain a sharp asymptotics for $n^{\KorElem_d}(\e)$ (it seems that this is impossible), but we have a logarithmic one.

\begin{Theorem}
Let $B_j$, $j\in\N$, satisfy  $\eqref{assum_param}$ with $c>0$, $\tau=1$, and  $0\leqslant r<\infty$. Then
\begin{eqnarray*}
\ln n^{\KorElem_d}(\e)=\bigl(1-(\e/\e_0)^2\bigr)\cdot \ln d +o(\ln d),\quad d\to\infty,\quad \text{for every}\quad \e\in(0,\e_0).
\end{eqnarray*}
\end{Theorem}
\textbf{Proof.}\quad The sequence $B_j$, $j\in\N$, satisfies Assumption $\ref{assum_tendstoinfty}$. Indeed,  since $(\beta_j)_{j\in\N}$ is  bounded, $(\zeta(\sigma_j))_{j\in\N}$ has a positive limit, and $\sum_{j=1}^d \beta_j\to\infty$, $d\to\infty$, the following sequence
\begin{eqnarray*}
	\dfrac{\max\limits_{j=1,\ldots,d}\tlambda^{B_j}_1}{\sum_{j=1}^d \tLambda^{B_j}}=\dfrac{\max\limits_{j=1,\ldots,d} \beta_j}{2\sum_{j=1}^d \beta_j \zeta(\sigma_j)},\quad d\in\N,
\end{eqnarray*}
tends to zero, and thus Assumption $\ref{assum_tendstoinfty}$ holds.

Let us define $\ell_j\defeq -\ln c+ \ln j$, $j\in\N$. The convergence
\begin{eqnarray*}
	\dfrac{1  }{\tLambda^{B_j}}\sum_{k=1}^\infty\tlambda^{B_j}_k\,\id\bigl(-\ln\tlambda^{B_j}_k-\ell_j\leqslant x\bigr)
	\to \id(x\geqslant 0),\quad j\to\infty,
\end{eqnarray*}
is established as just as in the  proof of Theorem \ref{th_nKorrElemde}.

Next, we consider the sequence from left-hand side of \eqref{th_nYd_UWconv_log_cond_W}, where we set $X_j=B_j$, $a_d\defeq 0$, $b_d\defeq \ln d$ $d\in\N$. Since $(\ell_j)_{j\in\N}$ increases, we set  $\signum\defeq 1$. For all $x\in\R$ and $d\in\N$ we have
\begin{eqnarray}\label{conc_tLambdaBj_conv}
	\dfrac{\sum_{j=1}^d\tLambda^{B_j} \id\bigl(\signum\cdot(\ell_j-a_d)\leqslant x b_d\bigr)}{\sum_{j=1}^d\tLambda^{B_j}}
	&=& \dfrac{\sum_{j=1}^d 2\beta_j \zeta(\sigma_j) \id\bigl(\ln j -\ln c\leqslant x \ln d\bigr)}{\sum_{j=1}^d2\beta_j \zeta(\sigma_j)}\nonumber\\
	&=& \dfrac{\sum_{j=1}^d \beta_j \zeta(\sigma_j) \id\bigl(j\leqslant c d^x\bigr)}{\sum_{j=1}^d \beta_j \zeta(\sigma_j)}.
\end{eqnarray}
It is easily seen that for $x<0$ and $x>1$ this fraction equals $0$ and $1$ respectively for all large enough $d$. Since, by the assumptions, $\sum_{j=1}^d \beta_j \zeta(\sigma_j)\to\infty$, $d\to\infty$, for $x=0$ the fractions in \eqref{conc_tLambdaBj_conv} tends to zero.  For $x\in(0,1]$ we can write
\begin{eqnarray*}
\sum_{j=1}^d \beta_j \zeta(\sigma_j) \id\bigl(j\leqslant c d^x\bigr)=\sum_{j=1}^{j_{d,x}} \beta_j \zeta(\sigma_j),\quad d\in\N,
\end{eqnarray*}
where $j_{d,x}\defeq\min\{d, \lfloor cd^x\rfloor\}$. Here $j_{d,x}\sim cd^x$ for $x\in(0,1)$, and $j_{d,x}\sim \min\{c,1\} d$ for $x=1$ as $d\to\infty$. On account of the assumptions on $\beta_j$, $\sigma_j$, $j\in\N$,
\begin{eqnarray}\label{conc_betajzetasjtauequals1}
\sum_{j=1}^d \beta_j \zeta(\sigma_j)\sim 	\sum_{j=1}^d \beta_j \sim \sum_{j=1}^d c/j \sim c \ln d, \quad d\to\infty.
\end{eqnarray}
From these remarks we conclude that
\begin{eqnarray*}
\dfrac{\sum_{j=1}^{j_{d,x}} \beta_j \zeta(\sigma_j)}{\sum_{j=1}^d \beta_j \zeta(\sigma_j)}\sim \dfrac{\ln j_{d,x}}{\ln d}\to x,\quad d\to\infty,\quad x\in(0,1].
\end{eqnarray*}
Thus the convergence \eqref{th_nYd_UWconv_cond_W} holds with $W(x)=\min\{x,1\}\id(x\geqslant 0)$, $x\in\R$.

All conditions of Theorem \ref{th_nYd_UWconv_log} hold. Substituting our parameters into \eqref{th_nYd_UWconv_log_conc},  we obtain
\begin{eqnarray*}
\ln n^{\KorElem_d}(\e)=  \ln\Bigl(\sum_{j=1}^d 2\beta_j\zeta(\sigma_j)\Bigr)+W^{-1}\bigl(1-(\e/\e_0)^2\bigr) \ln d+o(\ln d),\quad d\to\infty,
\end{eqnarray*}
for every $\e\in(0,\e_0)$ such that $1-(\e/\e_0)^2\in\ContSet(W^{-1})$. On account of \eqref{conc_betajzetasjtauequals1}, the first term from right-hand side equals  $\ln(\ln d)+\ln (2c)+ o(1)=o(\ln d)$, $d\to\infty$. Next, $W^{-1}(y)=y$, $y\in(0,1)$, and $\ContSet(W^{-1})=(0,1)$. Thus we obtain the required asymptotics 
\begin{eqnarray*}
\ln n^{\KorElem_d}(\e)=\bigl(1-(\e/\e_0)^2\bigr) \ln d+o(\ln d),\quad d\to\infty,
\end{eqnarray*}
for every $\e\in(0,\e_0)$.\quad $\Box$\\

Let us comment the behaviour of $n^{\KorElem_d}(\e)$ as $d\to\infty$. It is surprisingly that for the rather wide case $\tau<1$ the approximation complexity saves the linear growth as the function $d\mapsto C_\e d$, where $C_\e$ is a constant for fixed $\e$. At the same time only in the particular case $\tau=1$ growth of $n^{\KorElem_d}(\e)$ essentially differs. This becomes in fact sublinear and polynomial as $d\mapsto d^{q_\e}$, where $q_\e\in(0,1)$ is a constant for fixed $\e$.

\section*{Acknowlegments}
The work of the first named author was supported by DFG--SPbSU grant 6.65.37.2017, the RFBR grant 16-01-00258. 

Also the work of the first named author was also partially supported by the Government
of the Russian Federation (grant 074-U01), by Ministry of Science
and Education of the Russian Federation (GOSZADANIE 2014/190,
Projects No 14.Z50.31.0031 and No. 1.754.2014/K), by grant
MK-5001.2015.1 of the President of the Russian Federation.


\begin{thebibliography}{99}
\bibitem{ChenLi} X. Chen, W.V. Li, \textit{Small deviation estimates for some additive processes}, Proc. Conf. High Dimensional Probab. III, Progress in Probability, \textbf{55} (2003), Birkh\"{a}user, 225--238.

\bibitem{Hick} F. J. Hickernell, G. W. Wasilkowski, H. Wo\'zniakowski, \textit{Tractability of linear multivariate problems in the average-case setting}, in: A. Keller, S. Heinrich, H. Niederreiter (Eds.), Monte Carlo and Quasi-Monte Carlo Methods 2006, Springer, Berlin, 2008, pp. 461--493.

\bibitem{Khart} A. A. Khartov, \textit{Asymptotic analysis of average case approximation complexity of Hilbert space valued random elements}, J. Complexity, \textbf{31} (2015), 835--866.

\bibitem{LifPapWoz1} M. A. Lifshits, A. Papageorgiou, H. Wo\'zniakowski, \textit{Average case tractability of non-homogeneous tensor product problems}, J. Complexity, \textbf{28} (2012), 539--561.

\bibitem{LifZani1} M. A. Lifshits, M. Zani, \textit{Approximation complexity of additive random fields}, J. Complexity, \textbf{24} (2008), no. 3, 362--379.

\bibitem{LifZani2} M. A. Lifshits, M. Zani, \textit{Approximation of additive random fields based on standard information: Average case and probabilistic settings}, J. Complexity, \textbf{31} (2015), no. 5, 659--674.

\bibitem{NovWoz1}  E. Novak, H. Wo\'zniakowski, \textit{Tractability of Multivariate Problems. Volume I: Linear Information}, EMS Tracts Math. 6, EMS, Z\"urich, 2008.

\bibitem{Vaart} A. W. van der Vaart, \textit{Asymptotic Statistics}, Camb. Univ. Press, Cambridge, 1998.

\bibitem{WasWoz} G. W. Wasilkowski, H. Wo\'zniakowski, \textit{Average case optimal algorithms in Hilbert spaces}, J. Approx. Theory, \textbf{47} (1986), 17--25.

\bibitem{Zygmund} A. Zygmund, \textit{Trigonometric series}, Volume I, Camb. Univ. Press, Cambridge, 2002.
\end{thebibliography}
\end{document}